\documentclass[12pt]{amsart}

\usepackage{ifthen}

\usepackage{amssymb}

\usepackage{pstricks,pst-node,multido}

\makeatletter
\relax

\def\leadfill{\leaders\hbox to .7em{\hss.\hss}\hfill}

\renewcommand{\tocsection}[3]{%
  \indentlabel{\@ifnotempty{#2}{\ignorespaces#1 #2.\quad}}#3\mdseries{\leadfill}}

\def\@secnumfont{\bfseries}

\def\section{\@startsection{section}{1}%
   \z@{.7\linespacing\@plus\linespacing}{.5\linespacing}%
    {\normalfont\large\bfseries\boldmath\centering}}

\def\subsection{\@startsection{subsection}{2}%
  \z@{.5\linespacing\@plus.7\linespacing}{-.5em}%
  {\normalfont\bfseries\boldmath}}

\def\@captionheadfont{\sf} 
\def\@captionfont{\sf}

\def\c@tocdepth{1}

\makeatother

\DeclareMathOperator{\Hom}{Hom}
\DeclareMathOperator{\End}{End}
\DeclareMathOperator{\arcs}{arcs}
\DeclareMathOperator{\CG}{CG}

\newcommand{\cra}{\curvearrowright}
\newcommand{\equivp}{\equiv^+}

\newcommand{\un}{\backslash}
\newcommand{\PC}{P}
\newcommand{\GA}{GA}

\newcounter{daana} \newcounter{daanb}

\newcommand{\lista}[2]
{\begin{list}{#1}{%
\usecounter{daana}%
\setlength{\listparindent}{\parindent}%
\setlength{\leftmargin}{3.5em}%
\setlength{\itemindent}{0mm}%
\setlength{\labelsep}{.6em}%
\setlength{\labelwidth}{10em}%
\setlength{\itemsep}{0em}%
\setlength{\parsep}{1mm}%
\setlength{\listparindent}{\parindent}%
#2}}

\newcommand{\listb}[2]
{\begin{list}{#1}{%
\usecounter{daanb}%
\setlength{\rightmargin}{\leftmargin}%
\setlength{\listparindent}{\parindent}%
#2}}

\newcommand{\thspace}{2ex}

\newtheoremstyle{slanted}
  {\thspace}
  {\thspace}
  {\sl}
  {}
  {\bfseries}
  {.}
  {.5em}
  {}

\theoremstyle{slanted}

\newtheorem{theorem}{Theorem}[section]
\newtheorem{prop}[equation]{Proposition}

\newtheorem{lemma}[equation]{Lemma}
\newtheorem{coro}[equation]{Corollary}

\newtheoremstyle{roman}
  {\thspace}
  {\thspace}
  {\rm}
  {}
  {\it}
  {.}
  {.5em}
  {}

\theoremstyle{roman}

\newtheorem{lemdef}[equation]{Lemma/Definition}
\newtheorem{definition}[equation]{Definition}
\newtheorem{convention}[equation]{Convention}
\newtheorem{example}[equation]{Example}
\newtheorem{remark}[equation]{Remark}
\newtheorem{remarks}[equation]{Remarks}
\newtheorem{notation}[equation]{Notation}

\numberwithin{equation}{section}

\newcommand{\daanmath}{\mathbb}
\newcommand{\zz}{\daanmath{Z}}

\newcommand{\cc}{\daanmath{C}}

\newcommand\ra{\rightarrow}
 
\renewcommand{\>}{\rangle}

\newcommand{\barr}[1]{\overline{#1}}
\newcommand{\be}{\begin{equation}}
\newcommand{\ee}{\end{equation}}
\newcommand{\col}{\text{\upshape :\ }}

\begin{document} 
%

\title[Garside groupoid structures on the braid group]{A class of Garside groupoid structures on the pure braid group}

\author{Daan Krammer}
\address{Department of Mathematics, University of Warwick, Coventry CV4 7AL, UK}
\email{daan@maths.warwick.ac.uk}



\date{27 March 2006.}

\begin{abstract} We construct a class of Garside groupoid structures on the pure braid groups, one for each function (called labelling) from the punctures to the integers greater than~1. The object set of the groupoid is the set of ball decompositions of the punctured disk; the labels are the perimeters of the regions. Our construction generalises Garside's original Garside structure, but not the one by Birman--Ko--Lee. As a consequence, we generalise the Tamari lattice ordering on the set of vertices of the associahedron.
\end{abstract}

\maketitle


\begin{center} \parbox{.7\textwidth}{\tableofcontents} \end{center}

\section{Introduction}

In \cite{gar69} Garside solved the word and conjugacy problems in the braid groups. Birman, Ko and Lee \cite{bkl98} gave another solution, but it had something in common with Garside's approach. Indeed, both approaches were examples of what are now called {\em Garside structures\,} on the braid group. The general theory of Garside groups was developed mainly by Dehornoy; we mention three good overviews \cite{dehpar99}, \cite{deh00}, \cite{deh02}.

We propose a generalisation of the concept {\em Garside group\,} to {\em Garside groupoid}. It is not surprising, and the required proofs can easily be adapted from those in the literature covering the group case.

A short definition of Garside groupoids is as follows. A Garside groupoid is a group $G$ acting freely on the left on a lattice $L$ with the following properties.\begin{itemize}
\item The orbit set $G\backslash L$ is finite.
\item There exists an automorphism $\Delta$ (written on the right) of the lattice $L$ commuting with the $G$-action.
\item For any $x\in L$ the interval $[x,x\Delta]:=\{y\in L:x\leq y\leq x\Delta\}$ is finite.
\item The ordering on $L$ is generated by $x\leq y$ whenever $y\in[x,x\Delta]$.
\end{itemize} But this definition ought to be a theorem, and we shall use a different definition in the paper.

If the $G$-action on $L$ is transitive it is a Garside group. 

The main result of the paper is a class of new Garside (groupoid) structures on the pure braid group, one for each function, called labelling, from the set of punctures to~$\zz_{\geq 2}$. The objects for the involved groupoid are the ball decompositions of a punctured disk; the labels are the perimeters of the regions.

Garside's original Garside structure on the braid group \cite{gar69} is a special case of our construction, but Birman--Ko--Lee's structure \cite{bkl98} is not.

The set $[x,x\Delta]$ seems to be the vertex set of a polytope in a natural way, but we won't go into this. Among the polytopes obtained this way are the permutahedron (if all labels are~2) and the associahedron \cite{sta63}, \cite{lee89} (if all labels are~3).

As a consequence of our main result we get many finite lattices $[x,x\Delta]$. They seem new even in the particular case where all labels are $3$. In the more particular case where all triangles have a vertex in common, this was previously known by the work of Tamari \cite{tam51}, \cite{fritam67}, \cite[page~18]{gra78}. See subsection~\ref{ga87} for more details.

The paper is written in the language of Garside groups, not least because it allows us to use the results from \cite{deh00}, \cite{deh02} thus streamlining our proofs.

After the Garside structures by Garside and Birman--Ko--Lee there doesn't seem to be a need for any more of them. Our results are the byproducts of an investigation into surface mapping class groups. It would be interesting to know if mapping class groups are Garside, see~\cite{par05}.

In section~\ref{ga89} we construct a braid-like groupoid and a presentation for it. In section~\ref{ga90} we overview Garside groupoids in general. In section~\ref{ga91} we prove that the axioms for Garside groupoids are fulfilled by our braid groupoid. A short last section gives some more information about the finite lattices $[x,x\Delta]$ (there written $\Omega(xT)$) and especially the case where all labels are~3.

I would like to thank an anonymous referee who provided many suggestions for improvement.

\begin{remark} Depending on the labels, some of our proofs can be shortened, that is, some case-by-case proofs can deal with fewer cases. This will be immediately clear from the text. If all labels are $\geq 3$ then:
\begin{itemize}
\item The elementary relations (ER2) and (ER4) don't exist (figure~\ref{ga26}).
\item Cases~1 and 2 don't exist (figures~\ref{ga12} and \ref{ga16}).
\item Cases~(f), (g), (h) in the proof of lemma~\ref{ga70} (the cube condition) can be ignored.
\end{itemize}
If all labels are $3$ then:
\begin{itemize}
\item The elementary relations (ER2) and (ER4) don't exist (figure~\ref{ga26}).
\item Subsections~\ref{ga72} and \ref{ga93} can be ignored (subsection~\ref{ga94} replaces them).
\item Cases~(a), (f), (g), (h) in the proof of lemma~\ref{ga70} (the cube condition) can be ignored.
\end{itemize} \end{remark}

\section{A braid-like groupoid} \label{ga89}

\subsection{Introduction}

In this section we define a groupoid such that the automorphism groups of objects are finite index subgroups of the braid group. We give a presentation for it (generators and relations). It is a {\em complemented presentation\,} as required for Garside groupoids. The bulk of the section is there to prove that the presentation is correct.

\subsection{A braid groupoid}

Let $D$ be a disk, that is, a topological space homeomorphic to $\{z\in\cc:|z|\leq 1\}$. Let $P$ be a finite set of $n$ interior points of $D$ called {\em punctures}. Let $Q$ be a non-empty finite set of boundary points of $D$. Let $\ell$ be a map called {\em labelling\,} $\ell\col P\ra\zz_{\geq 2}$. For a reason which will become clear later we assume
\be |Q|-2=\sum_{p\in P}\big(\ell(p)-2\big). \label{ga1} \ee

Let $M$ be the {\em mapping class group\,} of
\[ (D,\mbox{orientation},P,Q,\ell). \]
In other words, $M=H/H_0$; here $H$ is the group of orientation preserving self-homeomorphisms $g$ of $D$ fixing $P$ setwise and the boundary $\partial D$ pointwise such that $\ell=\ell\circ g$; and $H_0$ is the component of identity in $H$.

Our definition of the {\em braid group\,} $B_n$ is precisely $M$ provided all labels are equal. In general, $M$ is isomorphic to the group of label preserving braids in $B_n$; it is a subgroup of $B_n$ of finite index
\[ n!\Big(\prod_{k\geq 2}n_k!\Big)^{-1} \]
where $n_k$ is the number of punctures of label~$k$.

\newcommand{\ninegon}{%
\pspicture[.4](-1.7,-1.7)(1.7,1.7)$ 
\psset{linewidth=.7pt, unit=.17mm, dimen=middle}
\SpecialCoor
\pscircle(0,0){100}
\pnode(100;  0){v0} \pnode(100; 40){v1}
\pnode(100; 80){v2} \pnode(100;120){v3}
\pnode(100;160){v4} \pnode(100;200){v5}
\pnode(100;240){v6} \pnode(100;280){v7}
\pnode(100;320){v8}
\pnode(100;60){e1} \pnode(100;100){e2} \pnode(100;220){e3} \pnode(100;260){e4} \pnode(100;300){e5} \pnode(100;340){e6}
\cnode*( 74;160){1.8pt}{p1} \uput[250]{0}(p1){3}
\cnode*( 74;320){1.8pt}{p1} \uput[230]{0}(p1){3}
\cnode*( 70; 68){1.8pt}{p1} \uput[180]{0}(p1){4}
\cnode*( 35; -120){1.8pt}{p3} \uput[-120]{0}(p3){5}
\pnode(0,0){c0} \pnode(135; 160){c1} \pnode(0,15){c2}
$\endpspicture}

\newcommand{\smallninegon}{%
\pspicture[.5](-.9,-.9)(.9,.9)$
\psset{linewidth=.7pt, unit=.09mm, dimen=middle, arcangle=45, ncurv=.8}
\SpecialCoor
\pscircle(0,0){100}
\pnode(100;  0){v0} \pnode(100; 40){v1}
\pnode(100; 80){v2} \pnode(100;120){v3}
\pnode(100;160){v4} \pnode(100;200){v5}
\pnode(100;240){v6} \pnode(100;280){v7}
\pnode(100;320){v8}
\cnode*( 74;160){1.6pt}{p1} 
\cnode*( 74;320){1.6pt}{p1} 
\cnode*( 70; 68){1.6pt}{p1} 
\cnode*( 35; -120){1.6pt}{p3} 
\pcarc(v3)(v5) \pcarc(v7)(v0)
$\endpspicture}

\newcommand{\ninedots}{%
\psline[linestyle=none, showpoints=true, dotsize=3.9pt, dotstyle=o](v0)(v1)(v2)(v3)(v4)(v5)(v6)(v7)(v8)}

\newcommand{\ninesmalldots}{%
\psline[linestyle=none, showpoints=true, dotsize=3.5pt, dotstyle=o](v0)(v1)(v2)(v3)(v4)(v5)(v6)(v7)(v8)}

\newcommand{\ninegonplus}{%
\ninegon \SpecialCoor \psset{arcangle=45, ncurv=.8, linewidth=.7pt, unit=.2mm}
\pcarc(v3)(v5) \pcarc(v7)(v0)
\psline[linestyle=none](v3)(v5)(v0)(v7)
}

\begin{figure}
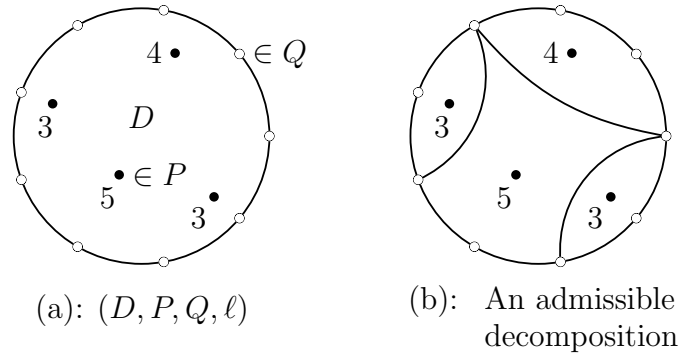

\[ \begin{array}{ccc} 
\SpecialCoor \ninegon \rput(c2){D} \uput[0]{0}(p3){\in P} \uput[0]{0}(v1){\in Q} \ninedots &&  
\ninegonplus \pcarc[arcangle=20](v0)(v3) \ninedots \\ 
\text{(a): $(D,P,Q,\ell)$} &\hspace{2em}&
\begin{tabular}{ll}\rule{0mm}{6mm} (b): & An admissible \\ & decomposition \end{tabular}
\end{array} \]
\caption{The punctured disk and an admissible decomposition.\label{ga4}}
\end{figure}

See figure~\ref{ga4}(a) for a picture of $(D,P,Q,\ell)$; each puncture $p$ is labelled $\ell(p)$.

By a {\em ball complex\,} $K$ we mean a CW complex whose cell attaching maps are injective. The cells of dimensions (respectively) 0,1,2 are called (respectively) {\em vertices, edges, regions}. By a {\em ball decomposition\,} of a topological space $X$ we mean a ball complex whose underlying topological space is $X$, and whose cells are subsets of $X$, and whose cell attaching maps are identity.

\begin{definition} An {\em admissible decomposition\,} is a ball decomposition of $D$ such that the following hold:\begin{itemize}
\item The vertex set is $Q$.
\item The 1-skeleton doesn't contain any punctures.
\item Each region (=\,2-dimensional cell) contains precisely one puncture $p$ and moreover, the label $\ell(p)$ equals the number of edges of that region.
\end{itemize}
\end{definition}

See figure~\ref{ga4}(b) for a picture of an admissible decomposition. For example, the puncture in the middle of the region with 4 edges is labelled~4.

\begin{remarks}  (a). The boundary of $D$ is in the $1$-skeleton of any admissible decomposition and consists of $|Q|$ edges which we shall call {\em boundary edges}. All other edges will be called {\em arcs}.

Every element of $Q$ is a vertex, even if there isn't any arc emanating from it.

(b). Recall our assumption (\ref{ga1}):
\be |Q|-2=\sum_{p\in P}\big(\ell(p)-2\big). \label{ga2} \ee
Equation (\ref{ga2}) is equivalent to saying that admissible decompositions exist. This follows because the number of even-dimensional cells minus the number of odd-dimensional cells is 1 (the Euler characteristic of the disk).
\end{remarks}

\begin{definition} By $L$ we will denote the set of isotopy classes relative $P$ (that is, $H_0$-orbits) of admissible decompositions.
\end{definition}

As is usual we shall abusively confuse admissible decompositions (or arcs) with their isotopy classes. Most of the time we are dealing with isotopy classes without saying we are.

The braid-like group $M$ acts on $L$, by convention on the left. It is clear that the action is {\em free\,} ($gx=x\Rightarrow g=1$ for all $g\in M$, $x\in L$) with a finite number of orbits.

Recall that a groupoid is a category all of whose morphisms are invertible. A group is then a groupoid with just one object. By $\Hom_G(A,B)$ or $G(A,B)$ we shall denote the set of morphisms in a category $G$ from an object $A$ to an object $B$. Its elements are said to have {\em source\,} object $A$ and {\em target\,} object $B$, and are said to go from $A$ to~$B$. The order of multiplication will be
\begin{align*} \Hom_G(A,B)\times\Hom_G(B,C) &\longrightarrow\Hom_G(A,C),  \\
(f,g) &\longmapsto fg. \end{align*}

\begin{definition} \label{ga82} We define a groupoid $G$ as follows. The object set is $G_0=M\backslash L$. 

We let $M$ act on cartesian powers $L^n$ diagonally: \[ g(x_1,\ldots,x_n)=(gx_1,\ldots,gx_n). \] An $M$-orbit in $L^n$ will be written $M(x_1,\ldots,x_n)$.

The hom-set $\Hom_G(A,B)$ is defined to be the set of $M$-orbits in $L\times L$ meeting (hence contained in) $A\times B$. Equivalently, $\Hom_G(Mx,My):=\{M(x,gy)\mid g\in M\}$.

The multiplication in $G$ is defined by $M(x,y)\circ M(y,z):=M(x,z)$. One readily checks that this is well-defined. Note the order of multiplication.

Note that this construction of $G$ could be applied whenever a group (in our case $M$) acts on a set (in our case $L$).
\end{definition}

\begin{remarks} \label{ga97} (a). It is easier to think of the group $M$ acting on the set $L$ rather than the groupoids constructed in definition~\ref{ga82}. Groupoids are just a language and don't add anything to our reasoning. It is a good idea not to lose sight of the group $M$ acting on the set~$L$.

(b). Clearly, for $x,y\in L$, one has $Mx=My\Longleftrightarrow x$ and $y$ are isotopic {\em not\,} relative $P$. In other words, the object set of $G$ can be identified with the set of isotopy classes, {\em not\,} relative $P$, of ball decompositions of $D$ as follows.
\begin{itemize}
\item The vertex set is $Q$.
\item For any $k\geq 2$, the number of $k$-gons equals the number of punctures with label~$k$.
\end{itemize}
\end{remarks}

\subsection{Elementary morphisms in $G$}

\begin{figure}
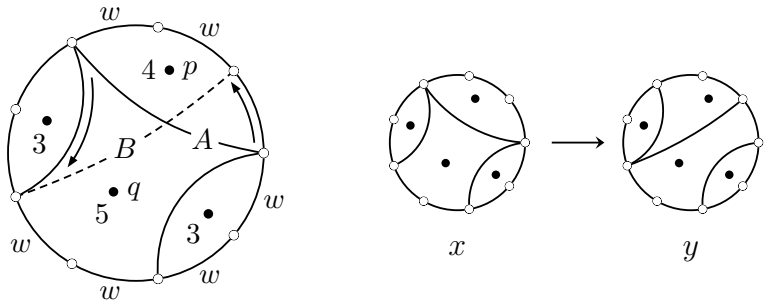

\begin{align*}
\SpecialCoor \small  \ninegonplus 
{\psset{labelsep=3pt} 
\nput{60}{e1}{w} \nput{100}{e2}{w} \nput{220}{e3}{w} 
\nput{260}{e4}{w} \nput{300}{e5}{w} \nput{340}{e6}{w} 
\nput{0}{p1}{p} \nput{0}{p3}{q}}
\pcarc[arcangle=20](v0)(v3) \ncput*[framesep=2pt, npos=.25]{\!A}
\pcarc[arcangle=10, linestyle=dashed, dash=3pt 2pt](v1)(v5)
\ncput*[framesep=2pt, npos=.53]{B}
\psset{arrowsize=1pt 3, arrowlength=1.4, arrowinset=0}
\psarc{<-}(c1){79.2}{-38}{8} \psarc{->}(c0){79}{5}{37} 
\ninedots &&
  \begin{array}{ccc}
  \smallninegon \SpecialCoor \pcarc[arcangle=28, linewidth=.7pt](v0)(v3) 
       \ninesmalldots & 
  \pspicture(0,0)(.6,0) \psset{arrowsize=2pt 4, 
        arrowlength=.7, linewidth=.7pt} 
        \psline{->}(0,0)(.7,0) \endpspicture &
  \smallninegon \SpecialCoor 
  \pcarc[arcangle=10, linewidth=.7pt](v1)(v5) \ninesmalldots \\ 
  \rule{0mm}{6mm}x && y 
  \end{array}
\end{align*}
\caption{An elementary pair $(x,y)$. See definition~\ref{ga98}.\label{ga5}}
\end{figure}

\begin{definition}[Elementary morphisms] \label{ga3} (See figure~\ref{ga5} and ignore the $w$'s for the moment). Let $A$ be an arc of an admissible decomposition~$x$. There is a unique ball decomposition $t$ of $D$ whose edges are the edges of $x$ except~$A$. Let $R$ denote the region of $t$ containing~$A$. Let $B$ be the arc obtained from $A$ by rotating both endpoints of $A$ in positive direction over one edge along the boundary of $R$, moving the rest of $A$ along in a continuous way such that it never meets any punctures. There is a unique admissible decomposition $y$ of $D$ whose edges are the edges of $t$ and~$B$. We call $(x,y)$ an {\em elementary pair\,} and $M(x,y)\in \Hom_G(Mx,My)$ an {\em elementary morphism}. In diagrams we draw elementary morphisms by solid arrows.
\end{definition}

For future reference we next define the diagonal morphisms, which generalise the elementary morphisms.

\begin{definition}[Diagonal morphisms] \label{ga98} Let $x$ be an admissible decomposition. Let $A_1,\ldots, A_k$ be $k$ distinct arcs of $x$. There is a unique ball decomposition $t$ of $D$ whose edges are the edges of $x$ except $A_1,\ldots, A_k$. For $1\leq i\leq k$ let $R_i$ denote the region of $t$ containing~$A_i$. Let $B_i$ be the arc obtained from $A_i$ by rotating both endpoints of $A_i$ in positive direction over one edge along the boundary of $R_i$, moving the rest of $A_i$ along in a continuous way such that it never meets any punctures. There is a unique admissible decomposition $y$ of $D$ whose edges are the edges of $t$ and $B_1,\ldots,B_k$. The morphism $M(x,y)\in \Hom_G(Mx,My)$ is called a {\em diagonal morphism}. If $k$ is the maximal value $n-1$ then we also call it a {\em full diagonal morphism}. We recover the elementary morphisms by taking $k=1$.
\end{definition}

For every $x\in L$ there are exactly $n-1=|P|-1$ elementary pairs $(x,y)$ (respectively, $(y,x)$), corresponding to the $n-1$ arcs of $x$ which are rotated in positive (respectively, negative) direction to obtain $y$.

In remark~\ref{ga97}(b) we saw that for the mere purpose of drawing an object of $G$, one can omit the punctures. What about morphisms?

It often happens that a hom-set $\Hom_G(Mx,My)$ contains at most one elementary morphism. This is at least the case if $Mx\neq My$. Moreover, if $Mx=My$ then the elementary morphisms from $Mx$ to $My$ are in bijection with those arcs of $x$ which lie on two bigons in $x$; rotating such an arc gives rise to the corresponding elementary morphism, see (\ref{ga100}). Thanks to this observation we shall often simplify our diagrams in $G$ by omitting the punctures, so that the arcs won't look like spaghetti.

\newcommand{\foura}{{%
\psset{unit=.4mm, linewidth=.7pt}
\raisebox{.7ex}{\pspicture[.5](-10,-13)(10,13)
\SpecialCoor \degrees[4]
\pnode(10;1){q0} \pnode(10;2){q1}
\pnode(10;3){q2} \pnode(10;4){q3}
\pnode(16;3.5){r}
\pscircle[dimen=middle](0,0){10}
\endpspicture}}}

\newcommand{\fourb}{{%
\psset{unit=.3mm, linewidth=.7pt}
\pspicture(-20,-20)(20,23)
\SpecialCoor \degrees[4]
\pnode(20;1){q0} \pnode(20;2){q1} \pnode(20;3){q2} \pnode(20;4){q3}
\pnode(0,0){c}
\pnode(-5, 14){p1} \pnode( 5, 14){p2} \pnode(-5,-14){p3} \pnode( 5,-14){p4}
\pscircle*(0,9){1.8pt} \pscircle*(0,-9){1.8pt}
\pscircle[dimen=middle](0,0){20}
\endpspicture}}

\newcommand{\drawdotsc}{%
\SpecialCoor \psline[linestyle=none, showpoints=true, dotstyle=o, dotsize=3.7pt](q0)(q1)(q2)(q3)}

\begin{example} \label{ga20}
Consider the case where the labels are $(3,3)$. Then there is a bijection $f\col\zz\ra L$. Here are some values of $f$.
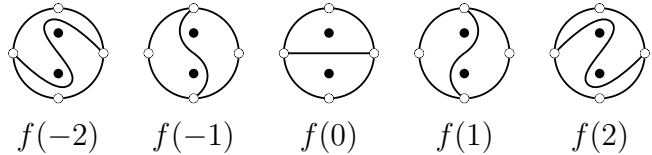
\begin{gather*} \psset{unit=1.2mm, linewidth=.7pt,  arrowsize=2pt 4, arrowlength=.7, ncurv=1, arcangle=55} \SpecialCoor \psmatrix[colsep=6mm, rowsep=0mm]
\fourb \pscurve(q3)(p1)(p4)(q1) \drawdotsc &
\fourb \pcarc(c)(q0) \pcarc(c)(q2) \drawdotsc &
\fourb \psline(q1)(q3) \drawdotsc &
\fourb \pcarc(q0)(c) \pcarc(q2)(c)  \drawdotsc &
\fourb  \pscurve(q3)(p3)(p2)(q1) \drawdotsc \\
f(-2) & f(-1) & f(0) & f(1) & f(2)
\endpsmatrix \label{ga108} \end{gather*}
We have $M\cong\zz$. A generator of $M$ takes $f(n)$ to $f(n+2)$. So $M\backslash L$, which is also $G_0$, the set of objects of $G$, has precisely two elements $x_0$, $x_1$ where $x_k=\{f(2n+k)\mid n\in\zz\}$. 

In remark~\ref{ga97}(b) we saw that elements of $G_0$ are given by pictures without punctures. Such pictures of $x_k$ are as follows.
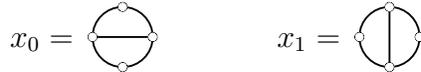
\[ \psset{unit=1.1mm, linewidth=.7pt,  arrowsize=2pt 4, arrowlength=.7, nodesep=6mm} \SpecialCoor
x_0={}\ \foura \psline(q1)(q3) \drawdotsc \hspace{4em}
x_1={}\ \foura \psline(q0)(q2) \drawdotsc \]

The elementary pairs are the pairs $(f(n),f(n+1))$. There are two elementary morphisms $g\col x_0\ra x_1$ and $h\col x_1\ra x_0$. Indeed $g$ is represented by the elementary pairs $(f(2n),f(2n+1))$, and $h$ by $(f(2n-1),f(2n))$. Note that $g$ and $h$ are not each other's inverses because $f(n)\neq f(n+2)$.

Consider $G(x_0,x_0)$, the set of morphisms in $G$ from $x_0$ to $x_0$. There is a bijection $r\col\zz\ra G(x_0,x_0)$ given by $r(n)=M(f(0),f(2n))$, that is, $r(n)$ is represented by the elementary pair $(f(0),f(2n))$. It restricts to a bijection $r\col\zz_{\geq0}\ra G^+(x_0,x_0)$.
\end{example}

\subsubsection*{A flash forward} We define a relation $\leq$ on $L$ by putting $x\leq y$ if and only if there is a sequence $x=x_0,x_1,\ldots,x_n=y$ such that $(x_i,x_{i+1})$ is an elementary pair for all $i$. We shall observe that $\leq$ is an ordering on $L$.  In example~\ref{ga20}, the ordering on $L$ is given by $f(k)\leq f(\ell)$ $\Leftrightarrow$ $k\leq\ell$. One of our results, corollary~\ref{ga67}, states that $(L,\leq)$ is a lattice.

\subsection{Elementary relations in $G$} \label{ga83}

\newcommand{\dota}{\pspicture(0,0) \pscircle*(0,0){1.6pt} \endpspicture}
\newcommand{\linea}{\psset{linewidth=.8pt, arrowsize=2pt 4, arrowlength=.7, dotstyle=o, dotsize=3.6pt, linestyle=solid}}
\newcommand{\dasha}{\linea \psset{linestyle=solid, linewidth=1.8pt, dash=2pt 1.5pt}}

\newcommand{\erone}{%
\pspicture[.5](0,0)(18mm,6mm) \linea \psset{yunit=.3mm, xunit=.75mm}
\rput( 5,17){\dota}  \rput( 5,3){\dota} 
\rput(19,17){\dota}  \rput(19,3){\dota}
\pnode( 0,20){v1}  \pnode(10,20){v2}
\pnode( 0, 0){v5}  \pnode(10, 0){v6}
\pnode(14,20){v3}  \pnode(24,20){v4}
\pnode(14, 0){v7}  \pnode(24, 0){v8}
\psline(v1)(v5) \psline(v2)(v6)
\psline(v4)(v8) \psline(v3)(v7)
\dasha \psset{arcangle=40}
\ncarc{v1}{v2} \ncarc{v3}{v4}
\ncarc{v6}{v5} \ncarc{v8}{v7}
\endpspicture}

\newcommand{\ertwo}{%
\pspicture(0,0)(0,0) \linea \psset{unit=.07mm} \SpecialCoor
\pnode(100; 27){v4}  \pnode(100; 63){v3}
\pnode(100;117){v2}  \pnode(100;153){v1}
\pnode(100;207){v5}  \pnode(100;243){v6}
\pnode(100;297){v7}  \pnode(100;333){v8}
\psarc(0,0){100}{(v4)}{(v3)} \psarc(0,0){100}{(v2)}{(v1)}
\psarc(0,0){100}{(v5)}{(v6)} \psarc(0,0){100}{(v7)}{(v8)}
\dasha
\psarc(0,0){100}{(v3)}{(v2)} \psarc(0,0){100}{(v1)}{(v5)}
\psarc(0,0){100}{(v6)}{(v7)} \psarc(0,0){100}{(v8)}{(v4)}
\rput(0,0){\dota} \rput(70,0){\dota} \rput(-70,0){\dota}
\endpspicture}

\newcommand{\ereightdots}{{\psset{showpoints=true, linestyle=none} \psline(v1)(v2)(v3)(v4)(v5)(v6)(v7)(v8)}}

\newcommand{\ersixdots}{{\psset{showpoints=true, linestyle=none} \psline(q0)(q1)(q2)(q3)(q4)(q5)}}

\newcommand{\erfourdots}{{\psset{showpoints=true, linestyle=none} \psline(q0)(q1)(q2)(q3)}}

\newcommand{\ersixgon}{%
{\pspicture(0,0)(0,0) \psset{unit=.06mm, dimen=middle} \degrees[120]
\pnode(100; 23){q0}  \pnode(100; 37){q1}
\pnode(100; 63){q2}  \pnode(100; 77){q3}
\pnode(100;103){q4}  \pnode(100; -3){q5} 
\psarc(0,0){100}{23}{37} \psarc(0,0){100}{63}{77} \psarc(0,0){100}{103}{-3} 
\dasha \psarc(0,0){100}{-3}{23} \psarc(0,0){100}{37}{63} \psarc(0,0){100}{77}{103}
\rput(65;10){\dota}  \rput(65;50){\dota}  \rput(55;90){\dota}
\endpspicture}}

\newcommand{\erfourgon}{%
{\pspicture(0,0)(0,0) \psset{unit=.07mm, dimen=middle}
\pnode(100; 60){q0} \pnode(100;120){q1}
\pnode(100;240){q2} \pnode(100;300){q3}
\pnode( 20, 20){r0} \pnode(-20, 20){r1}
\pnode(-20,-20){r2} \pnode( 20,-20){r3}
\psarc(0,0){100}{60}{120} \psarc(0,0){100}{240}{300} 
\dasha \psarc(0,0){100}{120}{240} \psarc(0,0){100}{300}{60}
\rput(0,0){\dota}  \rput(70,0){\dota}  \rput(-70,0){\dota}
\endpspicture}}

\begin{figure}
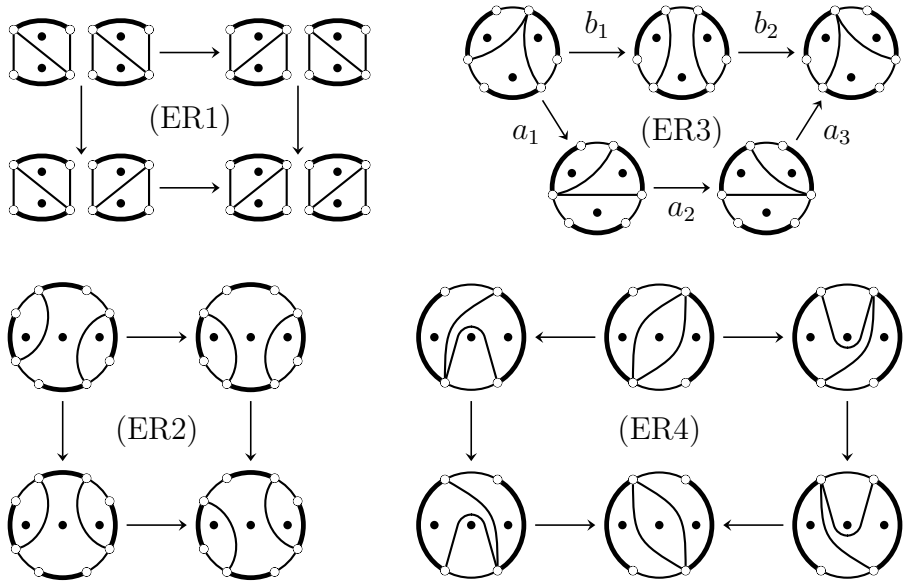
 \SpecialCoor \linea
\centering
\vspace{3mm} \psmatrix[rowsep=6mm] 
\psframebox[linestyle=none, framesep=0mm]{\pspicture[.5](-9mm,-4mm)(38mm,22mm)
\psset{unit=.9mm}
\rput(16, 10){\text{(ER1)}}
\rput(32,0){\rnode{n4}{\erone}} \psline(v2)(v5) \psline(v4)(v7) \ereightdots
\rput(0,0){\rnode{n3}{\erone}} \psline(v4)(v7) \psline(v1)(v6) \ereightdots
\rput(32,20){\rnode{n2}{\erone}} \psline(v2)(v5) \psline(v3)(v8) \ereightdots
\rput(0, 20){\rnode{n1}{\erone}} \psline(v1)(v6) \psline(v3)(v8) \ereightdots
\psset{nodesep=4pt, arrows=->, linewidth=.6pt}
\ncline{n1}{n2} \ncline{n1}{n3}
\ncline{n2}{n4} \ncline{n3}{n4}
\endpspicture}
\hspace{10mm}
\psframebox[framesep=6mm, linestyle=none]{\pspicture[.5](0mm,0mm)(45mm,18mm) \psset{unit=2.25mm}
\rput(10,3.3){\text{(ER3)}}
\pnode( 0,8){n1}  \pnode(10,8){n2}
\pnode(20,8){n3}  \pnode( 5,0){n4}
\pnode(15,0){n5}
\psset{arcangle=30}
\rput(n1){\ersixgon} \pcarc(q0)(q2) \pcarc(q4)(q0) \ersixdots
\rput(n2){\ersixgon} \pcarc(q4)(q0) \pcarc(q1)(q3) \ersixdots
\rput(n3){\ersixgon} \pcarc(q1)(q3) \pcarc(q5)(q1) \ersixdots
\rput(n4){\ersixgon} \pcarc(q0)(q2) \psline(q5)(q2) \ersixdots
\rput(n5){\ersixgon} \pcarc(q5)(q1) \psline(q2)(q5)  \ersixdots
\psset{nodesep=7.5mm, arrows=->, linewidth=.6pt}
$\ncline{n1}{n2}  \naput{b_1}
\ncline{n2}{n3}  \naput{b_2}
\ncline{n5}{n3}  \nbput{a_3}
\ncline{n1}{n4}  \nbput{a_1}
\ncline{n4}{n5}  \nbput{a_2}
$\endpspicture} \\
\psframebox[linestyle=none, framesep=7mm]{\pspicture(0mm,0mm)(25mm,25mm)
\psset{unit=1.25mm}
\rput(10, 10){\text{(ER2)}} 
\pnode( 0,20){n1} \pnode(20,20){n2}
\pnode( 0, 0){n3} \pnode(20, 0){n4}
\psset{arcangle=60, ncurv=.8}
\rput(n1){\ertwo} \pcarc(v2)(v5) \pcarc(v7)(v4) \ereightdots
\rput(n2){\ertwo}  \pcarc(v7)(v4) \pcarc(v1)(v6) \ereightdots
\rput(n3){\ertwo}  \pcarc(v2)(v5) \pcarc(v8)(v3) \ereightdots
\rput(n4){\ertwo}  \pcarc(v1)(v6) \pcarc(v8)(v3) \ereightdots
\psset{nodesep=8.5mm, arrows=->, linewidth=.6pt}
\ncline{n1}{n2} \ncline{n1}{n3}
\ncline{n2}{n4} \ncline{n3}{n4}
\endpspicture}
\hspace{12mm}
\psframebox[framesep=7mm, linestyle=none]{\pspicture(0mm,0mm)(50mm,25mm) \psset{unit=12.5mm}
\rput(2,1){\text{(ER4)}}
\pnode(0,2){n6} \pnode(2,2){n1}
\pnode(4,2){n2} \pnode(0,0){n5}
\pnode(2,0){n4} \pnode(4,0){n3}
\psset{arcangle=30, ncurv=1.1, unit=.07mm, linearc=27}
\rput(n1){\erfourgon} \ncarc{q0}{q2} \ncarc{q2}{q0}  \erfourdots
\rput(n2){\erfourgon} \ncarc{q0}{q2} \psline(q0)(r3)(r2)(q1) \erfourdots
\rput(n3){\erfourgon} \ncarc{q3}{q1} \psline(q0)(r3)(r2)(q1) \erfourdots
\rput(n6){\erfourgon} \ncarc{q2}{q0} \psline(q2)(r1)(r0)(q3)  \erfourdots
\rput(n5){\erfourgon} \ncarc{q1}{q3} \psline(q2)(r1)(r0)(q3)  \erfourdots
\rput(n4){\erfourgon} \ncarc{q1}{q3} \ncarc{q3}{q1} \erfourdots
\psset{nodesep=8.5mm, arrows=->, linewidth=.6pt}
\ncline{n1}{n2}  \ncline{n2}{n3}
\ncline{n3}{n4}  \ncline{n1}{n6}
\ncline{n6}{n5}  \ncline{n5}{n4}
\endpspicture} \endpsmatrix
\caption{Elementary relations. See definition~\ref{ga111}.\label{ga26}}
\end{figure}

\begin{definition} \label{ga111} An {\em elementary relation\,} is a relation in $G$ as displayed in figure~\ref{ga26}.

{\it Explanation.} (a). A fat line in an object stands for any nonnegative number of edges. If the number of edges is zero then the endpoints of the fat line are identical.

(b). In relations (ER2--4) any two objects of the relation differ only in a disk (a union of three regions) which is the only part depicted. In relation (ER1) two disks are involved, which may be disjoint or have a vertex or edge in common.

(c). All morphisms in elementary relations are elementary morphisms, and are therefore depicted by solid arrows.

With help of the punctures displayed in figure~\ref{ga26}, one observes that elementary relations really are relations in $G$, that is, they are commuting diagrams.
\end{definition}

\newcommand{\erfivespec}{%
\pspicture(0,0)(0,0) \SpecialCoor  \psset{unit=.06mm} \degrees[20] \pnode(100; 3){q1} \pnode(100; 7){q2} \pnode(100;11){q3} \pnode(100;15){q4} \pnode(100; 19){q0} \pscircle[dimen=middle](0,0){100} 
\psdots[dotsize=3.1pt](28;15)(70;1)(70;9)
\endpspicture}

\newcommand{\erfivedots}{{\psdots[dotstyle=o](q0)(q1)(q2)(q3)(q4)}}

\newcommand{\simplefive}{
\SpecialCoor \linea \psset{dotsize=3.4pt, dotstyle=*}
\psframebox[framesep=6mm, linestyle=none]{\pspicture[.5](0mm,0mm)(45mm, 18mm)  \psset{unit=2.25mm, arcangle=35, linewidth=.8pt}
\pnode( 0,8){n1}  \pnode(10,8){n2}
\pnode(20,8){n3}  \pnode( 5,0){n4}
\pnode(15,0){n5} 
\rput(n1){\erfivespec} \ncline{q3}{q1} \ncline{q1}{q4}  \erfivedots
\rput(n2){\erfivespec}\ncline{q1}{q4} \ncline{q4}{q2} \erfivedots
\rput(n3){\erfivespec}\ncline{q4}{q2}  \ncline{q2}{q0} \erfivedots
\rput(n4){\erfivespec}\ncarc{q3}{q0} \ncline{q3}{q1}  \erfivedots
\rput(n5){\erfivespec}\ncarc{q3}{q0} \ncline{q2}{q0} \erfivedots
\psset{arrows=->, nodesep=7mm, linewidth=.6pt}
\ncline{n1}{n2} \ncline{n2}{n3} \ncline{n5}{n3} \ncline{n1}{n4} \ncline{n4}{n5}
\endpspicture}}

\newcommand{\pent}{\pspicture(0,0)(0,0) \SpecialCoor  \psset{unit=3mm} \degrees[20] \pnode(1; 3){v1} \pnode(1; 7){v2} \pnode(1;11){v3} \pnode(1;15){v4} \pnode(1; 19){v5} \pspolygon[linewidth=.6pt](v1)(v2)(v3)(v4)(v5) \endpspicture}

\newcommand{\ppent}{%
\psframebox[linestyle=none, framesep=0pt]
{\pspicture[.4](-1.8,-.3)(1.8,1.5)
\psset{linewidth=.6pt, unit=1.5mm, arrowsize=2pt 4, arrowlength=.7, arrows=c-c}
\pnode( -5,0){n1}
\pnode(  5,0){n2}
\pnode(-10,8){n3}
\pnode(  0,8){n4}
\pnode( 10,8){n5}
\SpecialCoor
\rput(n1){\pent} \psline(v1)(v3) \psline(v3)(v5)
\rput(n2){\pent} \psline(v3)(v5) \psline(v5)(v2) 
\rput(n5){\pent} \psline(v5)(v2) \psline(v2)(v4) 
\rput(n4){\pent} \psline(v2)(v4) \psline(v4)(v1)
\rput(n3){\pent} \psline(v4)(v1) \psline(v1)(v3) 
\endpspicture}}

\begin{figure}
\begin{align*}
\begin{array}{cc} \\ \simplefive \\ \text{(a)} \end{array} && \psset{linewidth=.6pt, arrowsize=2pt 4, arrowlength=.7, nodesep=4mm}
\begin{array}{c@{}l}
\ppent \psset{arrows=->} \ncline{n3}{n4} \ncline{n4}{n5} \ncline{n3}{n1} \ncline{n1}{n2} \ncline{n2}{n5} & 
\text{(b): right} \\ \\
\ppent \psset{arrows=->} \ncline{n3}{n4} \ncline{n4}{n5} \ncline{n3}{n1} \ncline{n1}{n2} \ncline[arrows=<-]{n2}{n5} &
\text{(c): wrong}
\end{array}
\end{align*}
\caption{See example~\ref{ga112}.\label{ga109}}
\end{figure}

\begin{example} \label{ga112} (See figure~\ref{ga109}). If all labels are $3$ then the elementary relation (ER3) simplifies to part~(a) of figure~\ref{ga109}. On removing the punctures we obtain part~(b) of the figure. Diagram~(c) is obtained from diagram~(b) by reversing the right most arrow, and is therefore {\em not\,} a relation in $G$, certainly not an elementary relation. One toggles between right and wrong by either replacing each little $5$-gon by its mirror image, or by reversing all arrows.
\end{example}

\subsection{Presentations of groupoids}

There is a straightforward generalisation of the concept of group presentation (generators and relations) to presentations of categories and groupoids. We mention only that the object set is given in advance and is not affected by the presentation.

Recall that the braid group $B_n$ is the mapping class group of $(D,P,Q)$. It can also be given by the {\em Artin presentation\,} with generators $\sigma_i$ ($1\leq i<n$) and relations
\begin{align}
\text{$4$-gon:} && \sigma_{i}\, \sigma_{j} &=\sigma_{j}\, \sigma_{i} \qquad \text{whenever $|i-j|>1$,} \label{ga14} \\
\text{$6$-gon:} && \sigma_{i+1}\,\sigma_{i}\,\sigma_{i+1} &= \sigma_{i}\,\sigma_{i+1}\,\sigma_{i} \label{ga15}
\end{align} 
for all appropriate indices.

Our aim is to prove that $G$ is presented by elementary morphisms and elementary relations. Rather than building things from scratch we use the Artin presentation of the braid group as our starting point.

Let $G^*$ be the groupoid presented by all elementary morphisms and elementary relations. So we want to prove $G^*=G$.

\subsection{The case of only bigons} \label{ga17}

We now consider the special case where all labels are 2, so that every region of an admissible decomposition is a bigon. This makes a nice example, and will also be used later on.

\begin{lemma} \label{ga6} If all labels are~$2$, then $G^*=G$.
\end{lemma}

\newcommand{\bigons}{%
\psframebox[linestyle=none, framesep=0pt]{%
\pspicture[.5](-2.5,-.7)(2.5,.7) \SpecialCoor
\psset{linewidth=.6pt, unit=.4mm}
\cnode*(-50,0){2.1pt}{p1} \cnode*(-30,0){2.1pt}{p2}
\cnode*(-10,0){2.1pt}{p3} \cnode*( 10,0){2.1pt}{p4}
\cnode*( 30,0){2.1pt}{p5} \cnode*( 50,0){2.1pt}{p6}
\pnode(  0, 15){q0} 
\pnode( 0, -15){q1}
\pnode( 19 , 0  ) {h1}
\pnode( 9.5,-7.5) {h2}
\pnode(-9.5, 7.5) {h3}
\pnode(-19 , 0  ) {h4}
\pspolygon[linearc=14.9](-60,15)(-60,-15)(60,-15)(60,15)
\psline[linearc=7](q0)(40,6)(40,-6)(q1)
\psline[linearc=9](q0)(20,6)(20,-6)(q1)
\psline[linearc=7](q0)(-40,6)(-40,-6)(q1)
\psline[linearc=9](q0)(-20,6)(-20,-6)(q1)
\endpspicture}}

\begin{proof} Since all labels are equal, the automorphism group of any object of $G$ is the braid group $B_n$.

By our assumption, all labels are~$2$, so $G$ has only one object, and is equal to the braid group.

An elementary morphism looks like this:
\be
\bigons \SpecialCoor \psline[linewidth=.6pt,  dotsize=4.5pt]{o-o}(q0)(q1) \pspicture(0,0)(.7,0) \psset{linewidth=.6pt, arrowsize=2pt 4, arrowlength=.7} \psline{->}(0,0)(.7,0) \endpspicture
\bigons. \psset{linewidth=.6pt, unit=.4mm} \SpecialCoor \psline[linearc=5.9, dotsize=4.5pt]{o-o}(q0)(h1)(h2)(h3)(h4)(q1) \label{ga100} \ee
This is the pictorial interpretation of an Artin generator $\sigma_i$. Since $B_n$ is generated by the Artin generators, we conclude that $G$ is generated by the elementary morphisms.

One easily checks that the Artin $4$-gons (\ref{ga14}) are the elementary relations of the form (ER1) (figure~\ref{ga26}) and the Artin $6$-gons (\ref{ga15}) are the instances of (ER4). Relations (ER2) and (ER3) don't occur. The result follows.
\end{proof}

\subsection{Fans}

From now on, and throughout the section, we fix a base vertex $q_0\in Q$. 

\begin{definition} An admissible decomposition $x\in L$ is a {\em fan\,} if each arc of $x$ contains $q_0$. Its $M$-orbit $Mx\in G_0$ is also called a fan (confusion need not arise). Let $F$ denote the full subgroupoid of $G$ with the fans for objects ({\em full\,} means that every morphism of $G$ between fans is a morphism of~$F$).
\end{definition}

\begin{definition} By an {\em elementary fan morphism\,} we mean a morphism $M(x,y)\in\Hom_F(Mx,My)$ between two fans $Mx,My$ where an arc in $x$ is rotated in positive direction, not by one ``click'' as in elementary morphisms, but by the least positive number of clicks such that the result is again a fan (written $y$). In diagrams we denote elementary fan morphisms by dashed arrows. 
\end{definition}

\begin{convention} \label{ga30} From now on, with the exception of (\ref{ga31}), we will follow the following convention for depicting vertices in objects in diagrams. 

Of course: in a diagram, arrows go from objects to objects; in each object, there are vertices and edges.

In pictures of objects of $G$ (rather than the objects themselves) we don't speak of edges but rather {\em apparent edges}, which stand for a (usually unspecified) nonnegative number of edges. In other words, not all vertices may be depicted. If an apparent edge stands for zero edges, then its endpoints coincide. Apparent arcs always stand for a single arc.

A black dot $\bullet$ always stands for the base vertex $q_0$, but we allow the possibility that $q_0$ is depicted by a different symbol or not at all.

Note also that this convention will become more restricted from~\ref{ga73}.
\end{convention}

Let $D'$ denote the quotient of $D$ by contracting the union of all boundary edges not containing $q_0$, to a point. Define $P'$, $Q'$, $q_0'$ the obvious way, and let $\ell'\col P'\ra\zz_{\geq 2}$ be constant~$2$. Then $(D',P',Q',\ell')$ is again a collection of data as we started with, and indeed of the only-bigons form which we studied in subsection~\ref{ga17}. 

Define $G'$ accordingly. The edge contraction defines an isomorphism $F\ra G'$. But in lemma~\ref{ga6} we have seen a presentation for $G'$; on passing it through the isomorphism $F\cong G'$ we obtain a presentation for $F$. The generators in this presentation are precisely the elementary fan morphisms, and we shall call its relations the {\em fan relations}; they are just the Artin relations in disguise.

The $6$-gon fan relation is depicted in figure~\ref{ga24} using convention~\ref{ga30}. The labels $\sigma_i$ remind us that it is essentially an Artin relation. The perimeters of the regions are $k,\ell,m$.

\newcommand{\fan}{%
\pspicture(0,0)(0,0) \SpecialCoor \psset{unit=.065mm, linewidth=.7pt}
\pnode(100; 90){q0}  \pnode(100;240){q1}
\pnode(100;300){q2}  \pnode( 45;270){p2}
\pnode( 63;165){p1}  \pnode( 63; 15){p3}
\pscircle[dimen=middle](0,0){100}
\psset{showpoints=true, dotsize=4pt}
\psline[dotstyle=o](q1)(q0)(q2) \psline(q0)(q0)
\endpspicture}

\newcommand{\labs}[3]{%
{\everypsbox{}
\rput{0}(p1){#1} \rput{0}(p2){#2} \rput{0}(p3){#3}}}

\begin{figure}
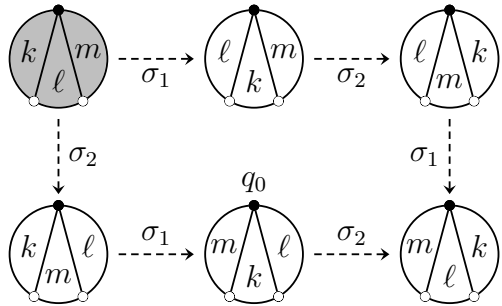
 
\[ \psframebox[linestyle=none, framesep=0pt]{\pspicture(-.7,-.7)(5.9,3.3) $ \SpecialCoor \psset{unit=26mm}
\pnode(0,1){n1}  \pnode(1,1){n2}
\pnode(2,1){n3}  \pnode(0,0){n4}
\pnode(1,0){n5}  \pnode(2,0){n6}
\pscircle[fillstyle=solid, fillcolor=lightgray, linestyle=none](n1){6.5mm}
\rput{0}(n1){\fan} \labs{k}{\ell}{m}
\rput{0}(n2){\fan} \labs{\ell}{k}{m}
\rput{0}(n3){\fan} \labs{\ell}{m}{k}
\rput{0}(n6){\fan} \labs{m}{\ell}{k}
\rput{0}(n5){\fan} \labs{m}{k}{\ell} \uput[90]{0}(q0){q_0}
\rput{0}(n4){\fan} \labs{k}{m}{\ell}
\psset{arrowsize=2pt 4, arrowlength=.7, linewidth=.7pt, nodesep=4mm, arrows=->, linestyle=dashed, dash=3pt 2pt, nodesep=8mm, labelsep=3pt} 
\ncline{n1}{n2} \nbput{\sigma_1}
\ncline{n2}{n3} \nbput{\sigma_2}
\ncline{n3}{n6} \nbput{\sigma_1}
\ncline{n1}{n4} \naput{\sigma_2}
\ncline{n4}{n5} \naput{\sigma_1}
\ncline{n5}{n6} \naput{\sigma_2}
$ \endpspicture} \]
\caption{The 6-gon fan relation.\label{ga24}}
\end{figure}

Summarising:

\begin{lemma} \label{ga7} The groupoid $F$ is presented by the elementary fan morphisms and the fan relations.\qed
\end{lemma}

\begin{lemma} \label{ga9} The groupoid $G$ is generated by the elementary morphisms.
\end{lemma}

\begin{proof} It is easy and left to the reader to show that one can walk from any admissible decomposition to a fan using elementary morphisms. In other words, every object of $G$ is isomorphic in $G^*$  to a fan. It remains to show that every morphism in $F$ is a product of elementary morphisms and their inverses (we say: it can be expressed in terms of elementary morphisms). Well, a morphism in $F$ can be expressed in terms of elementary fan morphisms (by lemma~\ref{ga7}) which in turn can be expressed in terms of elementary morphisms.\end{proof}

We now prepare for the proof that $G=G^*$. There is a natural surjective functor $G^*\ra G$. We don't yet know that it is injective, but we have at least a natural lifting of the elementary fan morphisms as follows. Every elementary fan morphism $u$ can be uniquely written as a product $u_1\cdots u_k$ of elementary morphisms where all arcs but one stay fixed. The expression $u_1\cdots u_k$ defines a morphism in $G^*$ which we write $a(u)$. The following is now clear.

\begin{lemma} \label{ga18} The statement $G=G^*$ is equivalent to saying that for every fan relation the following holds. On replacing each arrow $u$ in the fan relation by $a(u)$, one gets a relation in $G^*$ {\rm (}that is, a consequence of the elementary relations{\rm ).}\qed
\end{lemma}

\subsection{The case of only triangles} \label{ga94}

Before finishing the proof that $G=G^*$ in the general case we do the case where all labels are~$3$, which is easier yet shows all steps.

\begin{prop} \label{ga10} If all labels are~$3$ then $G^*=G$.
\end{prop}

\begin{proof} We prove the necessary and sufficient condition of lemma~\ref{ga18}. 

A 4-gon fan relation can be `tessellated' by four squares of the form (ER1); the details are easy and left to the reader.

Consider a 6-gon fan relation (see figure~\ref{ga24}). It necessarily looks like the dashed arrows in the following diagram.
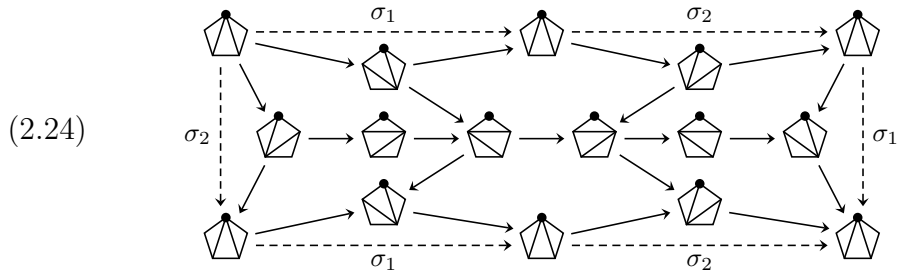
\newcommand{\penta}{\pspicture(0,0)(0,0) \SpecialCoor  \psset{unit=3mm} \degrees[20] \cnode*(1; 5){1.8pt}{v1} \pnode(1; 9){v2} \pnode(1;13){v3} \pnode(1;17){v4} \pnode(1; 1){v5} \pspolygon[linewidth=.6pt](v1)(v2)(v3)(v4)(v5) \endpspicture}
\newcommand{\pentaa}{\penta\psset{linewidth=.6pt, linearc=.1pt} \SpecialCoor \psline{c-c}(v3)(v1) \psline{c-c}(v1)(v4)}
\newcommand{\pentab}{\penta\psset{linewidth=.6pt, linearc=.1pt} \SpecialCoor \psline{c-c}(v5)(v3) \psline{c-c}(v3)(v1)}
\newcommand{\pentac}{\penta\psset{linewidth=.6pt, linearc=.1pt} \SpecialCoor \psline{c-c}(v2)(v4) \psline{c-c}(v4)(v1)}
\newcommand{\pentad}{\penta\psset{linewidth=.6pt, linearc=.1pt} \SpecialCoor \psline{c-c}(v4)(v2) \psline{c-c}(v2)(v5)}
\newcommand{\pentae}{\penta\psset{linewidth=.6pt, linearc=.1pt} \SpecialCoor \psline{c-c}(v2)(v5) \psline{c-c}(v5)(v3)}
%
\begin{gather} \text{\footnotesize \psframebox[linestyle=none, framesep=0pt]{%
\pspicture[.5](-48mm,-18mm)(48mm,18mm) \psset{yunit=.45mm, xunit=.7mm, arrowsize=2pt 4, arrowlength=.7, linewidth=.6pt, nodesep=4mm, arrows=->, labelsep=3pt} \SpecialCoor 
\pnode(-60,30){n1} \pnode(-30,20){n2} \pnode(0,30){n3} \pnode(30,20){n4} \pnode(60,30){n5} \pnode(-50,0){n6} \pnode(-30,0){n7} \pnode(-10,0){n8} \pnode(10,0){n9} \pnode(30,0){n10} \pnode(50,0){n11} \pnode(-60,-30){n12} \pnode(-30,-20){n13} \pnode(0,-30){n14} \pnode(30,-20){n15} \pnode(60,-30){n16}
\rput{0}(n1){\pentaa}\rput{0}(n3){\pentaa} \rput{0}(n5){\pentaa} \rput{0}(n12){\pentaa} \rput{0}(n14){\pentaa} \rput{0}(n16){\pentaa} \rput{0}(n4){\pentab} \rput{0}(n15){\pentab} \rput{0}(n6){\pentab} \rput{0}(n2){\pentac} \rput{0}(n11){\pentac} \rput{0}(n13){\pentac} \rput{0}(n8){\pentad} \rput{0}(n10){\pentad} \rput{0}(n7){\pentae} \rput{0}(n9){\pentae}
\ncline{n1}{n2} \ncline[offsetB=-2pt]{n2}{n3} \ncline{n3}{n4} \ncline[offsetB=-2pt]{n4}{n5} \ncline{n1}{n6} \ncline{n2}{n8} \ncline{n4}{n9} \ncline{n5}{n11} \ncline{n6}{n7} \ncline{n7}{n8} \ncline{n8}{n9} \ncline{n9}{n10} \ncline{n10}{n11} \ncline{n6}{n12} \ncline{n8}{n13} \ncline{n9}{n15} \ncline{n11}{n16} \ncline{n12}{n13} \ncline[offsetB=2pt]{n13}{n14} \ncline{n14}{n15} \ncline[offsetB=2pt]{n15}{n16}
\psset{linestyle=dashed, dash=3pt 2pt, offset=2pt} \ncline{n1}{n3} \naput{$\sigma_1$} 
\ncline{n3}{n5}  \naput{$\sigma_2$}
\ncline[nodesepB=4.5mm]{n5}{n16}  \naput{$\sigma_1$}
\psset{offset=-2pt}
\ncline[nodesepB=4.5mm]{n1}{n12}  \nbput{$\sigma_2$}
\ncline{n12}{n14}  \nbput{$\sigma_1$}
\ncline{n14}{n16} \nbput{$\sigma_2$}
\endpspicture}} \label{ga31} \end{gather} 
The triangles in this diagram show how to replace each elementary fan morphism $u$ of the fan relation by $a(u)$ consisting of two solid arrows. It is therefore our aim to prove that the twelve outer solid arrows form a relation in $G^*$. Well they do, because the six $5$-gons are elementary relations as one can check. In particular there are no wrong $5$-gons as in figure~\ref{ga109}(c). This settles the case of 6-gon fan relations.\end{proof}

\subsection{Complementary relations} \label{ga72}

\begin{convention} In addition to convention~\ref{ga30} the present subsection~\ref{ga72}  uses the following notation. If a vertex in an object is depicted by any symbol other than a white dot (namely, $\bigtriangleup, \bullet, 1,2,3,\ldots$) this indicates that the vertex is being repeated in nearby objects (same place, same symbol). Knowing this helps understand the diagrams. White dots $\circ$ denote other vertices which we decide to depict.
\end{convention}

\begin{definition} \label{ga113} Assume $n=3$, so that every admissible decomposition has $n-1=2$ arcs. We define two more special sorts of morphisms in $G$.

(a). Suppose the two arcs of $x\in L$ have a common vertex. Let $y\in L$ be the result of rotating both arcs of $x$ in $D$ over the same number ($\in\zz$) of clicks. The morphism $M(x,y)\in\Hom_G(Mx,My)$ is called a {\em central morphism\,} and notated by a solid arrow with a twiddle~$\stackrel{\sim}{\longrightarrow}$. Equivalently, a central morphism is any product of full diagonal morphisms, which were defined in \ref{ga3}, and their inverses. See figure~\ref{ga86} for pictures of central morphisms. 

(b). Let $a,b,c\in Q$ be three distinct vertices in this cyclic order counterclockwise. Let $x\in L$ be such that one of its arcs $A$ connects $a,b$ and the other arc $B$ connects $b,c$. Let $y\in L$ be obtained from $x$ by rotating $B$ counterclockwise by the least positive number of clicks such that it contains $a$. The morphism $M(x,y)\in\Hom_G(Mx,My)$ is called a {\em parabolic morphism}. Like elementary morphisms, they will be denoted by solid arrows; some parabolic morphisms are indeed elementary morphisms. The two non-twiddled arrows in (CR2) in figure~\ref{ga86} are examples of parabolic morphisms.

(c). The collective name for central and parabolic morphisms is {\em complementary morphisms}.
\newcommand{\radiusa}{4.8mm}
\newcommand{\radiusb}{6.8mm} 
\newcommand{\sizea}{20mm}
\newcommand{\sizeb}{40mm} 
\newcommand{\smallbox}[1]{\psovalbox[framesep=1.2pt, fillstyle=solid, linestyle=none, linewidth=.1pt]{\hspace*{-.3pt} #1\hspace*{-.3pt}}}
\newcommand{\fivea}{%
\pspicture(0,0)(0,0) \SpecialCoor \psset{unit=\radiusa} \degrees[20]
\pnode(1;5){q0}  \pnode(1;9){q1}
\pnode(1;13){q2} \pnode(1;17){q3}
\pnode(1;21){q4} \pnode(0,0){c} 
\pnode(1;3){q0a} \pnode(1;7){q1a} 
\pscircle[dimen=middle]{1} \endpspicture}
\newcommand{\squaredots}{%
{\SpecialCoor \psdots[dotstyle=triangle](q0)(q1)(q2)(q3)(q4)}}
\newcommand{\whitedots}{%
{\SpecialCoor \psdots(q0)(q1)(q2)(q3)(q4)}}
\newcommand{\colourdots}{%
{\SpecialCoor \psdots(q2)(q4) \psdots[dotstyle=triangle](q3)(q3)}}
\begin{figure}
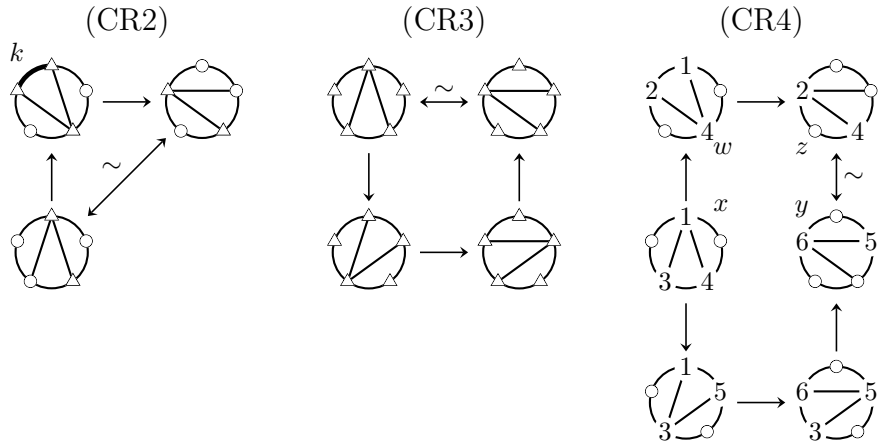

\begin{align*}
\psset{unit=\sizea, arrowsize=2pt 4, arrowlength=.7, nodesep=\radiusb, dotsize=5pt, dotstyle=o}
\setlength{\arraycolsep}{0mm}
\begin{array}{c@{\hspace{12mm}}c@{\hspace{12mm}}c}
\text{(CR2)} & \text{(CR3)} & \text{(CR4)} \\[3mm]
\psframebox[linestyle=none, framesep=\radiusa]{%
\pspicture[1](0mm,\sizea)(\sizea,\sizeb) \footnotesize $ \SpecialCoor
\pnode(0,1){n2} \pnode(0,2){n3} \pnode(1,2){n6}
\rput{0}(n2){\fivea} \psline(q2)(q0)(q3) \colourdots 
\psdots[dotstyle=triangle](q0) \psdots(q1) 
\rput{0}(n3){\fivea} \degrees[20] \psarc[linewidth=1.8pt](c){\radiusa}{5}{9} \uput[7]{0}(q1a){k} \psline(q0)(q3)(q1) \colourdots
\psdots[dotstyle=triangle](q0)(q1) 
\rput{0}(n6){\fivea} \psline(q3)(q1)(q4) \colourdots 
\psdots[dotstyle=triangle](q1) \psdots(q0) 
\psset{arrows=->, labelsep=2pt, linewidth=.6pt}
\ncline{n2}{n3}  \ncline{n3}{n6}
\ncline{<->}{n2}{n6} \naput{\sim}
$\endpspicture} &
\psframebox[linestyle=none, framesep=\radiusa]{%
\pspicture[1](0mm,0mm)(\sizea,\sizea) \footnotesize $ \SpecialCoor \psset{labelsep=2pt}
\pnode(0,0){n1}  \pnode(0,1){n2}
\pnode(1,0){n4}  \pnode(1,1){n5}  
\rput{0}(n1){\fivea} \psline(q4)(q2)(q0) \squaredots
\rput{0}(n2){\fivea} \psline(q2)(q0)(q3) \squaredots
\rput{0}(n4){\fivea} \psline(q1)(q4)(q2) \squaredots
\rput{0}(n5){\fivea} \psline(q3)(q1)(q4) \squaredots
\psset{arrows=->, linewidth=.6pt}
\ncline{n1}{n4}  \ncline{n2}{n1}
\ncline{n4}{n5}
\ncline{<->}{n2}{n5} \naput{\sim}
$\endpspicture} &
\psframebox[linestyle=none, framesep=\radiusa]{%
\pspicture[1](0mm,0mm)(\sizea,\sizeb) \footnotesize $ \SpecialCoor \degrees[20]
\pnode(0,0){n1}  \pnode(0,1){n2}
\pnode(0,2){n3}  \pnode(1,0){n4}
\pnode(1,1){n5}  \pnode(1,2){n6}
\rput{0}(n1){\fivea} \psline(q4)(q2)(q0) \whitedots
\rput(q2){\smallbox{3}} \rput(q4){\smallbox{5}} \rput(q0){\smallbox{1}}
\rput{0}(n2){\fivea} \psline(q2)(q0)(q3) \whitedots
\rput(q0){\smallbox{1}} \rput(q2){\smallbox{3}} \rput(q3){\smallbox{4}}
\uput[3]{0}(q0a){x}
\rput{0}(n3){\fivea} \psline(q0)(q3)(q1) \whitedots
\rput(q3){\smallbox{4}} \rput(q0){\smallbox{1}} \rput(q1){\smallbox{2}}
\uput[17]{0}(q3){w}
\rput{0}(n4){\fivea} \psline(q1)(q4)(q2) \whitedots
\rput(q4){\smallbox{5}} \rput(q1){\smallbox{6}} \rput(q2){\smallbox{3}}
\rput{0}(n5){\fivea} \psline(q3)(q1)(q4) \whitedots
\rput(q1){\smallbox{6}} \rput(q4){\smallbox{5}}
\uput[7]{0}(q1a){\smash{y}}
\rput{0}(n6){\fivea} \psline(q3)(q1)(q4) \whitedots
\rput(q3){\smallbox{4}} \rput(q1){\smallbox{2}}
\uput[13]{0}(q2){z}
\psset{arrows=->, labelsep=2pt, linewidth=.6pt}
\ncline{n1}{n4}  \ncline{n2}{n1}
\ncline{n2}{n3}  \ncline{n3}{n6}
\ncline{n4}{n5}
\ncline{<->}{n6}{n5} \naput{\sim}
$\endpspicture}
\end{array} \end{align*} 
\caption{Complementary relations. See definition~\ref{ga27}. \label{ga86}}
\end{figure}%
\end{definition}
\begin{definition} \label{ga27} Assume $n=3$. The {\em complementary relation\,} (CR1) reads $uv=w$ for any central morphisms $u,v,w$ (provided it is meaningful and true in $G$). 

The complementary relations (CR2--4) are defined in figure~\ref{ga86}. The fat apparent edge in (CR2) stands for $k\geq1$ edges. 

Finally, (CR5) are the relations $u=u_1\cdots u_k$ where $u$ is parabolic and $u_i$ is elementary, provided the relation is meaningful and true in $G$ and one arc is fixed along the path $u_1\cdots u_k$ (recall that there are only two arcs).
\end{definition}

\begin{example} Here are some examples of reading diagrams by the rules we have introduced.

(a). The non-twiddled arrows in (CR2) are not necessarily elementary. 

(b). The non-twiddled arrows in (CR3) are elementary. Three of the apparent boundary edges stand for single edges.

(c). Consider the parabolic arrow $x\ra w$ in (CR4). It is possible that the vertices are distinct and ordered (unexpectedly) $(1,3,2,4)$ in this cyclic order. The vertex $2$ is not necessarily repeated in $x$ nor $3$ in~$w$.

(d). Source and target of the twiddled arrow $z\ra y$ in (CR4) are not necessarily equal. They may involve different sets of depicted vertices, and the pictures of the disk may differ by a rotation.
\end{example}

\begin{lemma} \label{ga23} Assume $n=3$. The groupoid $G^*$ is presented by elementary and complementary morphisms and elementary and complementary relations.
\end{lemma}

\begin{proof} By definition, $G^*$ is presented by elementary morphisms and relations. For each complementary morphism $u$ we shall define a morphism $b(u)$ of $G^*$ (with the same source and target objects as $u$). The result will then be proved by showing that replacing each complementary arrow $u$ by $b(u)$ in the complementary relations yields relations in $G^*$ (that is, consequences of the elementary relations).

The twiddled arrow in (CR2) with $k=1$ is a central morphism $u_0$ which for the moment we call a simple morphism; we define $b(u_0)$ and $b(u_0^{-1})$ by that diagram. Any other central morphism $u$ is a product $u_1\cdots u_\ell$ of simple morphisms and their inverses and we define $b(u):=b(u_1)\cdots b(u_\ell)$.

A parabolic morphism $u$ can uniquely be written as a product $u_1\cdots u_k$ of elementary morphisms where all arcs but one stay fixed, as in (CR5); we define $b(u):=u_1\cdots u_k$.

Finally we prove that replacing $u\mapsto b(u)$ in complementary relations yields relations in $G^*$. For (CR1) and (CR5) this is trivial. For (CR2) this is an easy induction on $k$ and left to the reader. For (CR3) this follows from the $5$-gon (ER3) and (CR2) with $k=1$. In order to prove (CR4), cut it into three pieces along the diagonals $xy$ and $xz$ and apply (CR2), (CR1) and (CR3) to the pieces (respectively, from top to bottom).\end{proof}

\begin{remark} \label{ga58} Note that (CR4) generalises (ER3). In order for (CR4) to be true we had to insert a central morphism into a $5$-gon making it into a $6$-gon. Inserting the central morphism can be done anywhere in the $5$-gon (not necessarily the spot where (CR4) does) and we shall tacitly do so when necessary.
\end{remark}

\subsection{Presentation of $G$} \label{ga93}

We will use a variation of one of our lemmas as follows. Assume $n=3$. Let $u$ be an elementary fan morphism. We define a morphism $c(u)$ in $G^*$ as follows.\begin{itemize}
\item If $u$ is an elementary morphism we put $c(u)=u$.
\item Otherwise, in $G^*$ one can write $u$ uniquely as $u_1u_2$ (with one arc fixed throughout) where $u_1$ or $u_2$ is an elementary morphism and the other is a parabolic morphism. We define $c(u)=u_1u_2$.
\end{itemize}
The following mild variation of lemma~\ref{ga18} is now clear.

\begin{lemma} \label{ga29} Assume $n=3$. The statement $G=G^*$ is equivalent to the following.\begin{itemize}
\item For every $4$-gon fan relation, the condition of lemma~{\rm \ref{ga18}} holds.
\item For every $6$-gon fan relation, the following holds. On replacing each arrow $u$ in the fan relation by $c(u)$, one gets a relation in $G^*$ {\rm (}that is, a consequence of the elementary and complementary relations{\rm ).}\qed
\end{itemize}
\end{lemma}


\newcommand{\fourgon}{%
\pspicture(0,0)(0,0)
\psset{unit=4mm}
\SpecialCoor
\degrees[4]
\pnode(1;1){q0}  \pnode(1;2){q1}
\pnode(1;3){q2}  \pnode(1;4){q3}
\pnode(0,0){c}
\pscircle[dimen=middle, linewidth=.8pt](0,0){1}
\endpspicture}

\newcommand{\drawdotsb}{%
\SpecialCoor
\psline[linestyle=none, showpoints=true, dotsize=4pt, dotstyle=o](q1)(q2)(q3)
\pscircle*(q0){2pt}}


\newcommand{\dlinea}[2]{%
\SpecialCoor \degrees[20]
\psarc[doubleline=false, doublesep=1pt, linewidth=1.8pt](c){4mm}{#1}{#2}}

\newcommand{\fivegona}{%
\pspicture(0,0)(0,0)
\psset{unit=4mm}
\SpecialCoor
\degrees[20]
\pnode(1; 5){q0}  \pnode(1;9){q1}
\pnode(1;13){q2}  \pnode(1;17){q3}
\pnode(1;21){q4} 
\pnode(0,0){c}
\pscircle[dimen=middle, linewidth=.8pt](0,0){1}
\endpspicture}

\newcommand{\drawdotsa}{%
\SpecialCoor
\psline[linestyle=none, showpoints=true, dotsize=4pt, dotstyle=o](q1)(q2)(q3)(q4)
\pscircle*(q0){2pt}}


\newcommand{\dline}[2]{%
\SpecialCoor \degrees[28]
\psarc[doubleline=false, doublesep=1pt, linewidth=1.8pt](c){4.5mm}{#1}{#2}}

\newcommand{\sevengon}{%
\pspicture(-0,-0)(0,0)
\psset{unit=4.5mm}
\SpecialCoor
\degrees[28]
\pnode(1; 7){q0}  \pnode(1;11){q1}
\pnode(1;15){q2}  \pnode(1;19){q3}
\pnode(1;23){q4}  \pnode(1;27){q5}
\pnode(1; 3){q6}
\pnode(0,0){c}
\pscircle[dimen=middle, linewidth=.8pt](0,0){1}
\endpspicture}

\newcommand{\drawdots}{%
\SpecialCoor
\psline[linestyle=none, showpoints=true, dotsize=3.8pt, dotstyle=o](q1)(q2)(q3)(q4)(q5)(q6)
\pscircle*(q0){1.9pt}}

\begin{figure}[b]
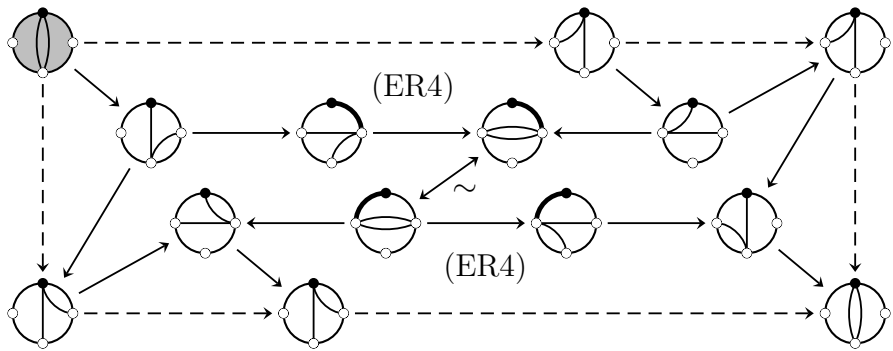

\centering
\vspace*{3ex} 
\psframebox[linestyle=none]{%
\pspicture(-60mm,-22mm)(60mm,22mm)
\psset{xunit=1.2mm, yunit=1mm, linewidth=.7pt}
\SpecialCoor
\pnode( -45, 18){v1}  \pnode(  15, 18){v2}
\pnode(  45, 18){v3}  \pnode( -33,  6){v4}
\pnode( -13,  6){v5}  \pnode(   7,  6){v6}
\pnode(  27,  6){v7}  \pnode(  45,-18){v14}
\pnode( -15,-18){v13} \pnode( -45,-18){v12}
\pnode(  33, -6){v11} \pnode(  13, -6){v10}
\pnode(  -7, -6){v9}  \pnode( -27, -6){v8}
\pscircle[linestyle=none, fillstyle=solid, fillcolor=lightgray](v1){4mm}
\rput(-4,12){\text{(ER4)}}
\rput(4,-12){\text{(ER4)}}
\psset{arcangle=45}
\rput{0}(v2){\fourgon}  \psline(q0)(q2) \pcarc(q0)(q1) \drawdotsb
\rput{0}(v3){\fourgon}  \psline(q0)(q2) \pcarc(q0)(q1) \drawdotsb
\rput{0}(v4){\fourgon} \psline(q0)(q2) \pcarc(q2)(q3) \drawdotsb
\rput{0}(v5){\fourgon} \psline(q1)(q3) \pcarc(q2)(q3) \dlinea{0}{5} \drawdotsb
\rput{0}(v7){\fourgon} \psline(q1)(q3) \pcarc(q0)(q1) \drawdotsb
\rput{0}(v8){\fourgon} \psline(q1)(q3) \pcarc(q3)(q0) \drawdotsb
\rput{0}(v10){\fourgon} \psline(q1)(q3) \pcarc(q1)(q2) \dlinea{5}{10} \drawdotsb
\rput{0}(v11){\fourgon} \psline(q0)(q2) \pcarc(q1)(q2) \drawdotsb
\rput{0}(v12){\fourgon}  \psline(q0)(q2) \pcarc(q3)(q0) \drawdotsb
\rput{0}(v13){\fourgon}  \psline(q0)(q2) \pcarc(q3)(q0) \drawdotsb
\degrees[360] \psset{arcangle=24} 
\rput{0}(v1){\fourgon} \pcarc(q0)(q2) \pcarc(q2)(q0) \drawdotsb
\rput{0}(v6){\fourgon} \pcarc(q1)(q3) \pcarc(q3)(q1) \dlinea{0}{5} \drawdotsb
\rput{0}(v9){\fourgon} \pcarc(q1)(q3) \pcarc(q3)(q1) \dlinea{5}{10} \drawdotsb
\rput{0}(v14){\fourgon} \pcarc(q0)(q2) \pcarc(q2)(q0) \drawdotsb
\psset{arrowsize=2pt 4, arrowlength=.7, nodesep=5.5mm, arrows=->, linewidth=.7pt}
\ncline{v1}{v4}   \ncline{v2}{v7} 
\ncline{v3}{v11}  \ncline{v4}{v5} 
\ncline{v4}{v12}  \ncline{v5}{v6} 
\ncline{v7}{v6}   \ncline{v7}{v3} 
\ncline{v8}{v13}  \ncline{v9}{v8} 
\ncline{v9}{v10}  \ncline{v10}{v11} 
\ncline{v11}{v14} 
\ncline{v12}{v8} 
\psset{labelsep=2pt}
\ncline{<->}{v6}{v9}   \naput[npos=.47]{$\sim$}
\psset{linestyle=dashed}
\ncline{v1}{v2}  \ncline{v2}{v3}   \ncline{v3}{v14}
\ncline{v1}{v12} \ncline{v12}{v13} \ncline{v13}{v14}
\endpspicture}
%
\caption{\label{ga12}Case 1: $\ell=2$.}
\end{figure}

\begin{figure}
\vspace{3ex} Case 2: $\ell>2$, $m=2$ \\[1ex]
\psframebox[linestyle=none]{%
\pspicture(-4.5,-1.75)(4.5,1.75)
\psset{unit=.95mm, linewidth=.7pt}
\SpecialCoor
\pnode( -40, 15){v1}  \pnode(   0, 15){v2}
\pnode(  40, 15){v3}  \pnode( -23,  0){v4}
\pnode(   0,  0){v5}  \pnode(  23,  0){v6}
\pnode( -40,-15){v7}  \pnode(   0,-15){v8}
\pnode(  40,-15){v9}
\pscircle[linestyle=none, fillstyle=solid, fillcolor=lightgray](v1){4mm}
\psset{arcangle=36}
\rput{0}(v1){\fivegona} \psline(q0)(q2) \pcarc(q4)(q0) \drawdotsa
\rput{0}(v3){\fivegona} \psline(q0)(q2) \pcarc(q0)(q1) \drawdotsa
\rput{0}(v4){\fivegona} \psline(q2)(q4) \pcarc(q4)(q0) \drawdotsa
\rput{0}(v5){\fivegona} \psline(q2)(q4) \pcarc(q1)(q2) \dlinea{5}{9} \drawdotsa
\rput{0}(v6){\fivegona} \psline(q1)(q3) \pcarc(q0)(q1) 
\dlinea{21}{5} \drawdotsa
\rput{0}(v7){\fivegona} \psline(q0)(q3)  \pcarc(q4)(q0) \drawdotsa
\rput{0}(v9){\fivegona} \psline(q0)(q3) \pcarc(q0)(q1) \drawdotsa
\degrees[360] \psset{arcangle=24} 
\rput{0}(v2){\fivegona} \pcarc(q0)(q2) \pcarc(q2)(q0) \drawdotsa 
\rput{0}(v8){\fivegona} \pcarc(q3)(q0) \pcarc(q0)(q3) \drawdotsa 
\psset{arrowsize=2pt 4, arrowlength=.7, nodesep=5.5mm, arrows=->, linewidth=.7pt}
\ncline{v1}{v4} \ncline[offsetB=0pt]{v4}{v7}
\ncline{v4}{v5}
\ncline{v3}{v6} \ncline[offsetB=0pt]{v6}{v9}
\psset{labelsep=2pt}
\ncline{<->}{v6}{v5}   \nbput[npos=.47]{$\sim$}
\psset{linestyle=dashed}
\ncline{v3}{v9} \ncline{v1}{v7} \ncline{v7}{v8}
\ncline{v8}{v9} \ncline{v1}{v2} \ncline{v2}{v3}
\endpspicture}
\\[3ex] Case 3: $\ell>2$, $m>2$ \\[1ex]
\psframebox[linestyle=none]{%
\pspicture(-4.8,-6.7)(4.8,6.7)
\psset{yunit=.75mm, xunit=1.05mm, linewidth=.8pt}
\SpecialCoor
\pnode( 40, 84){v1}
\pnode(  0, 84){v2}
\pnode(-40, 84){v3}
\pnode( 20, 68){v4}
\pnode(-20, 68){v5}
\pnode( 20, 40){v6}
\pnode(  0, 20){v7}
\pnode(-20, 40){v8}
\pnode( 34, -20){v9}
\pnode( 20,  0){v10}
\pnode(-20,  0){v11}
\pnode(-34, 20){v12}
\pnode( 20, -40){v13}
\pnode(  0, -20){v14}
\pnode(-20, -40){v15}
\pnode( 20, -68){v16}
\pnode(-20, -68){v17}
\pnode( 40, -84){v18}
\pnode(  0, -84){v19}
\pnode(-40, -84){v20}
\pscircle[linestyle=none, fillstyle=solid, fillcolor=lightgray](v2){4.5mm}
\rput{0}(v1){\sevengon} \psline(q2)(q0)(q4) \drawdots
\rput{0}(v2){\sevengon} \psline(q2)(q0)(q4) \drawdots
\rput{0}(v3){\sevengon} \psline(q3)(q0)(q4) \drawdots
\rput{0}(v4){\sevengon} \psline(q0)(q2)(q5) \drawdots
\rput{0}(v5){\sevengon} \psline(q0)(q4)(q2) \dline{7}{11} \drawdots
\rput{0}(v6){\sevengon} \psline(q1)(q5)(q2) \drawdots
\rput{0}(v7){\sevengon} \psline(q5)(q1)(q3) \dline{3}{7} \drawdots
\rput{0}(v8){\sevengon} \psline(q6)(q3)(q1) \dline{3}{7} \drawdots
\rput{0}(v9){\sevengon} \psline(q0)(q4)(q2) \drawdots
\rput{0}(v10){\sevengon} \psline(q6)(q2)(q4) \dline{7}{11} \drawdots
\rput{0}(v11){\sevengon} \psline(q1)(q5)(q3) \dline{3}{7} \drawdots
\rput{0}(v12){\sevengon} \psline(q0)(q3)(q5) \drawdots
\rput{0}(v13){\sevengon} \psline(q1)(q4)(q6) \dline{7}{11} \drawdots
\rput{0}(v14){\sevengon} \psline(q2)(q6)(q4) \dline{7}{11} \drawdots
\rput{0}(v15){\sevengon} \psline(q6)(q2)(q5) \drawdots
\rput{0}(v16){\sevengon} \psline(q0)(q3)(q5) \dline{3}{7} \drawdots
\rput{0}(v17){\sevengon} \psline(q0)(q5)(q2) \drawdots
\rput{0}(v18){\sevengon} \psline(q3)(q0)(q4) \drawdots
\rput{0}(v19){\sevengon} \psline(q5)(q0)(q3) \drawdots
\rput{0}(v20){\sevengon} \psline(q5)(q0)(q3) \drawdots
\psset{arrowsize=2pt 4, arrowlength=.7, nodesep=6mm, arrows=->, linewidth=.7pt}
\ncline{v1}{v9}
\ncline{v2}{v4}
\ncline{v2}{v5}
\ncline{v3}{v12}
\ncline{v4}{v1}
\ncline{v4}{v6}
\ncline{v5}{v3}
\ncline{v6}{v7}
\ncline{v7}{v11}
\ncline{v8}{v7}
\ncline[offsetB=-4pt]{v9}{v18}
\ncline{v10}{v9}
\ncline{v10}{v14}
\ncline{v12}{v11}
\ncline[offsetB=4pt]{v12}{v20}
\ncline{v14}{v13}
\ncline{v14}{v15}
\ncline{v15}{v17}
\ncline{v16}{v19}
\ncline{v17}{v19}
\ncline{v18}{v16}
\ncline{v20}{v17}
\psset{labelsep=2pt}
\ncline{<->}{v8}{v5}   \nbput[npos=.43]{$\sim$}
\ncline{<->}{v7}{v10}  \naput{$\sim$}
\ncline{<->}{v11}{v14} \naput{$\sim$}
\ncline{<->}{v16}{v13} \nbput[npos=.43]{$\sim$}
\psset{linestyle=dashed}
\ncline[offset= 4pt]{v1}{v18}
\ncline[offset=-4pt]{v3}{v20}
\ncline{v2}{v1}    \ncline{<-}{v3}{v2}
\ncline{v18}{v19}  \ncline{v20}{v19}
\psset{boxsize=9mm, nodesep=8.5mm, linearc=8.4mm, linewidth=.5pt, dash=3pt 2pt}
\ncbox{v5}{v8}    \ncbox{v7}{v10}
\ncbox{v11}{v14}  \ncbox{v13}{v16}
\endpspicture}
%
\caption{\label{ga16}}
\end{figure}

\begin{prop} \label{ga28} $G=G^*$.
\end{prop}

\begin{proof} As in the proof of proposition~\ref{ga10}, we verify the necessary and sufficient condition of lemma~\ref{ga29}.

First consider a $4$-gon fan relation $uv=wx$ where $u,v,w,x$ are elementary fan morphisms. Then there are $k,\ell\geq 1$ such that both $a(u),a(x)$ are a product of $k$ elementary morphism while $a(v),a(w)$ consist of $\ell$ of them. Then the fan relation $uv=wx$ can be tessellated by $k\ell$ $2$-cells of the form (ER1).

We turn to the $6$-gon fan relation (see figure~\ref{ga24}). It entirely takes place in a subdisk of $3$ regions so without loss of generality we may assume $n=3$. Then the definition of complementary morphisms and relations applies. Let $k,\ell,m\geq 2$ be the labels in the order defined by figure~\ref{ga24}. We consider three cases:
\begin{align*}
& \text{Case 1:}\quad \ell = 2,  \\
& \text{Case 2:}\quad \ell >2,\ m=2, \\
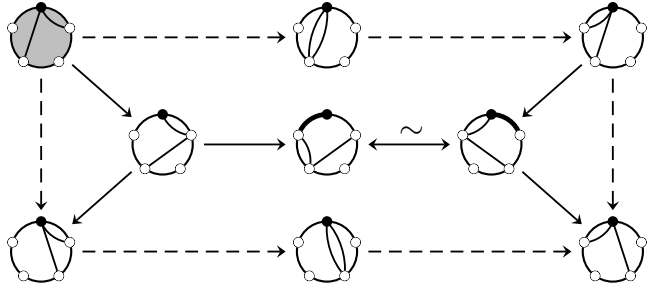
& \text{Case 3:}\quad \ell >2,\ m>2.
\end{align*}
The proofs of cases~1, 2, 3 can be found in (respectively) figures~\ref{ga12}, \ref{ga16},~\ref{ga16}. For convenience we have shaded one object of every fan $6$-gon.

In order to be able to read the diagrams, we need to know which vertices are repeated from object to object. The rule, which is different from the previous subsection, is as follows.

If the value of a label is prescribed (namely, $\ell=2$ in case~1 and $m=2$ in case~2) then each apparent edge of a corresponding region stands for one edge.

Any other label only has a determined lower bound. These lower bounds are $k\geq 2$ (all cases), $\ell\geq 3$ (cases 2 and~3), $m\geq 2$ (case~1), $m\geq 3$ (case~3). Consider one of these lower bounds, say ``label${}\geq r$''. Every corresponding region will have $r+1$ apparent edges, $r$ of which stand for one edge, the remaining {\em exceptional\,} one standing for any nonnegative number of edges. In most cases, it is clear or it doesn't matter which of the apparent edges is the exceptional one and we will not mark it to keep the diagrams simple. The exception is when both boundary edges containing $q_0$ are contained in the same region%
\footnote{It is rather flexible which vertex is depicted by a white dot, but not so $q_0$.}.
In that case and some other cases the exceptional apparent edge is drawn fat.

Consequently, the diagrams can almost be read as if every apparent edge stands for just one edge.

The main work is in checking that all bounded $2$-cells in our diagrams are relations in $G^*$, so that we may conclude that so are the unbounded $2$-cells. This cannot be done in print, but we mention the following. The triangles show how to express an elementary fan morphism $u$ (a dashed arrow) as $c(u)$ (solid arrows), a product of an elementary and a parabolic morphism in some order. Two $2$-cells in figure~\ref{ga12} are (ER4) as indicated. The $4$-gon in the middle of figure~\ref{ga16} (second diagram) is not an elementary or complementary relation but is easily shown to be a consequence of such. All remaining (bounded) $2$-cells are (CR4) but they may have the central morphism in a different place as explained in \ref{ga58}. This finishes the proof.\end{proof}

Note that contraction of the central arrows in the second diagram of figure~\ref{ga16} (they are enclosed in dashed ovals) recovers diagram~(\ref{ga31}), at least when taken as planar directed graph.

%
\section{Garside groupoids} \label{ga90}
%

\subsection{Introduction}

In this section we give a summary about Garside groupoids in general. It aims to be readable by people unfamiliar with Garside groups.

The literature is restricted to Garside groups rather than groupoids but it is easy to generalise the assertions and proofs of all results about Garside groups to groupoids. If you want to learn about Garside groupoids, the easiest thing to do is to read about Garside groups first, because their notation is simpler. Conceptually there is hardly any difference. For this reason, we begin by summarising Garside groups. Subsequently we give a summary of Garside groupoids, which of course is almost identical to the group case.

Good overviews on Garside groups are \cite{dehpar99}, \cite{deh00}, \cite{deh02}. Our summary is mainly based on \cite{deh00} but our terminology is mainly taken from \cite{deh02}.

\subsection{Garside groups} \label{ga102}

\begin{figure}[h]
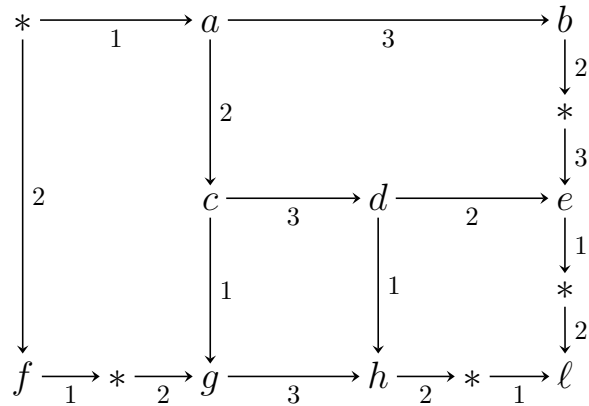
 \large \psset{nodesep=3pt, linewidth=.6pt, labelsep=3pt, arrowsize=2pt 4, arrowlength=.7, arrows=->}
\centering
\[
\psmatrix[colsep=4.5ex, rowsep=2.7ex]
{*} &{}& a &{}&{}&{}& b \\
&{}&{}&{}&{}&{}&{*} \\
&{}& c &{}& d &{}& e \\
{}&{}&{}&{}&{}&{}&{*} \\
f &{*}& g &{}& h &{*}& \ell 
\everypsbox{\scriptstyle}
\ncline{1,1}{1,3}  \nbput{1}
\ncline{1,3}{1,7}  \nbput{3}
\ncline{3,3}{3,5}  \nbput{3}
\ncline{3,5}{3,7}  \nbput{2}
\ncline{5,1}{5,2}  \nbput{1}
\ncline{5,2}{5,3}  \nbput{2}
\ncline{5,3}{5,5}  \nbput{3}
\ncline{5,5}{5,6}  \nbput{2}
\ncline{5,6}{5,7}  \nbput{1}
\ncline{1,1}{5,1}  \naput{2}
\ncline{1,3}{3,3}  \naput{2}
\ncline{3,3}{5,3}  \naput{1}
\ncline{3,5}{5,5}  \naput{1}
\ncline{1,7}{2,7}  \naput{2}
\ncline{2,7}{3,7}  \naput{3}
\ncline{3,7}{4,7}  \naput{1}
\ncline{4,7}{5,7}  \naput{2}
\endpsmatrix \]
\caption{\label{br104} Word reversing. See~example~\ref{ga105}.}
\end{figure}

\begin{example} \label{ga105} As an introduction to Garside groups we given an example. It is called the word reversing process and is at the basis of the theory of Garside groups.

Consider the braid group $B_n$ presented by generators $\sigma_i$ and relations (\ref{ga14}), (\ref{ga15}). The braid monoid $B_n^+$ is the submonoid of $B_n$ generated by the $\sigma_i$. Define an ordering $\leq$ on $B_n$ by $x\leq z$ if and only if $z=xy$ for some $y\in B_n^+$. For $x,y\in B_n$, let $x\vee y$ denote the {\em join\,} of $x,y$, that is, the least common upper bound, provided it exists. It is known that
 the joins $\sigma_i\vee\sigma_j$ exist and are given by
\be
\begin{aligned} 
\sigma_i\vee\sigma_j &=\sigma_i \, \sigma_j \, \sigma_i \qquad \qquad && |i-j|=1 \\
\sigma_i\vee\sigma_j &=\sigma_i \, \sigma_j &&|i-j|>1.
\end{aligned}
\label{breq37} \ee
We write $i$ for $\sigma_i$ and aim to compute $2\vee 13$ (and to prove that this join exists). See figure~\ref{br104}. By (\ref{breq37}) we have
\[ g=a\vee f, \quad e=c\vee b, \quad h=g\vee d, \quad \ell=h\vee e. \]
(In particular, these joins exist.) Going backwards through these identities, we find
\begin{align*} \ell &=h\vee e = g\vee d\vee e = g\vee e = g \vee c \vee b = g\vee b = a\vee f\vee b = f\vee b. \end{align*}
We conclude $2\vee 13=132312$. In algebraic notation the same process looks as follows.
\begin{align*} 
\underline{ 2^{-1} \, 1 } \, 3  &\ra 1 \, 2 \, 1^{-1} \, \underline{ 2^{-1} \, 3 }  \ra 1 \, 2\, \underline{ 1^{-1} \, 3 } \, 2\, 3^{-1} \, 2^{-1} \\
&\ra   1 \, 2\, 3 \, \underline{ 1^{-1} \, 2 } \, 3^{-1} \, 2^{-1}  \ra  1 \, 2\, 3 \, 2\, 1\, 2^{-1} \, 1^{-1} \, 3^{-1} \, 2^{-1}
\end{align*}
Thus, word reversing aims to remove all occurrences $a^{-1}b$ where $a,b$ are (positive) letters.
\end{example}

\subsubsection{Complemented presentations} Let $A$ be a finite set. Let $A^{-1}$ be a disjoint copy of $A$ and let $a\mapsto a^{-1}$ denote a bijection from $A$ to $A^{-1}$, whose inverse is also written $a\mapsto a^{-1}$. Let $W^+$ denote the free monoid on $A$ and $W$ the free monoid on $A\cup A^{-1}$. We call elements of $A$ {\em letters}, elements of $W$ {\em words\,} and elements of $W^+$ {\em positive words}. If $u\in W$ is a word, its inverse $u^{-1}$ is defined by inverting all letters in $u$ and reversing their order. Note that $u^{-1}u\neq 1$ unless $u=1$.

A {\em complement\,} on $A$ is a map $f$ which takes any pair of distinct $a,b\in A$ to a positive word $f(a,b)\in W^+$. We also call $(A,f)$ a complement.

The monoid $G^+$ (respectively, group $G$) associated with the above complement is defined as the monoid (respectively, group) presented by generators $A$ and relations
\be a\,f(a,b)=b\,f(b,a) \label{ga106} \ee
for all distinct $a,b\in A$.

We have a natural map $W^+\ra G^+$. We write $u\equivp v$ (where $u,v\in W^+$) if $u,v$ have the same image in $G^+$.

\begin{example} We shall show that the Artin presentation of the braid group is a complemented presentation.

In the above notation, assume that $A=\{\sigma_1,\ \ldots,\ \sigma_{n-1}\}$ and $f(\sigma_i,\sigma_j)=\sigma_j$ if $|i-j|>1$ and $f(\sigma_i,\sigma_j)=\sigma_j\,\sigma_i$ if $|i-j|=1$. Then the relations (\ref{ga106}) are precisely the Artin relations (\ref{ga14}), (\ref{ga15}) as promised.
\end{example}

\subsubsection{The cube condition}

For letters $a,b$ and words $u,v$ write
\[ u \, a^{-1} \, b \, v \ra u f(a,b) \, f(b,a)^{-1} v. \]
This is what we called word reversing in example~\ref{ga105}. For words $u,v$ write $u\cra v$ if there are $n\geq 0$ and $u_i$ such that
\[ u=u_0\ra u_1\ra\cdots\ra u_n=v. \]

Suppose $u,v,u_0,v_0$ are positive words such that $u^{-1}v\cra v_0u_0^{-1}$. It can be shown that $v_0$ depends only on $u,v$; it will be written $u\un v$ (``$u$ under $v$''). It may happen though that $u\un v$ is not defined, if $u_0,v_0$ with the required properties don't exist. Note $a\un b=f(a,b)$ for all $a,b\in A$.

We say that our complement satisfies the {\em cube condition}%
\footnote{In \cite{deh02} this is called the {\em weak cube condition on letters}.} 
if
\be (a\un b)\un (a\un c)\equivp (b\un a)\un (b\un c) \label{ga104} \ee
for all $a,b,c\in A$; here, and henceforth, we take an equation like (\ref{ga104}) to mean that either both sides are not defined, or both sides are defined and equivalent as asserted.

If you find the cube condition mysterious then you may like a more appealing version of it in terms of characteristic graphs, which is explained in remark~\ref{ga96}.

\subsubsection{Main theorem on Garside groups}

\begin{definition} A monoid is called {\em atomic\,} if for each element $x$ there is a natural number $n$ such that $x$ cannot be written as a product of $n$ or more nontrivial elements.
\end{definition}

\begin{definition} An {\em automorphism\,} of a complement $(A,f)$ (written on the  right) is an automorphism $\phi$ of the associated monoid $W^+$ such that $A\phi=A$ and $(f(a,b))\phi=f(a\phi,b\phi)$ for all $a,b\in A$.
\end{definition}

\begin{definition} Retain the above notation. A {\em Garside element\,} is an element $\Delta\in G^+$ satisfying the following.
\lista{{\rm(GE\arabic{daana})}}{\setlength{\leftmargin}{4em}}
\item For all $a\in A$ we have $\Delta\in a\,G^+$.
\item There exists an automorphism $\phi$ of $(A,f)$ such that for all $a\in A$ we have $a\Delta=\Delta(a\phi)$ (in $G^+$).
\end{list}
If $\Delta=uv$ (in $G^+$) then we call $u$ a {\em simple element}.
\end{definition}

\begin{definition} A {\em Garside group tuple\,} consists of a complement $(A,f)$ and a Garside element $\Delta$ such that $G^+$ is atomic, and the cube condition (\ref{ga104}) is satisfied.

In this case, we call $G$ a {\em Garside group\,} and $G^+$ a {\em Garside monoid\,} but technically one works with Garside group tuples.
\end{definition}

The following theorem is a rather arbitrary collection of the many good properties Garside groups have, to give the reader an idea. For proofs we refer to \cite{dehpar99}, \cite{deh00}, \cite{deh02} and the references cited therein.


\begin{theorem} \label{ga101} Let $(A,f,\Delta)$ be a Garside group tuple. Then the following hold.%
\lista{{\rm(G\arabic{daana})}}{\setlength{\leftmargin}{3em}}
\item Let $[u]$ denote the image in $G^+$ of a morphism $u$ in $W^+$. If $u,v\in W^+$ then $[u\un v]$ depends only on $[u]$ and $[v]$ and will therefore be written $[u]\un[v]$. 
\item Let us put an ordering on $G^+$ by $a\leq c$ if and only if $c=ab$ for some $b\in G^+$. Then $G^+$ is a lattice. That is, any two elements have a join {\rm (}a least common upper bound{\,\rm )} and a meet {\rm (}a greatest common lower bound{\,\rm ).} The join $u\vee v$ of $\{u,v\}\subset G^+$ is $u(u\un v)=v(v\un u)$.
\item The natural map $G^+\ra G$ is injective.
\item Every element of $G$ is a quotient of two elements in the image of the map $G^+\ra G$.
\item If $\Delta=uv$ in $G^+$ then $v$ is a simple element.
\item {\rm (}Greedy form{\rm ).} Every element $u\in G^+$ can uniquely be written $u=u_1\cdots u_n$ $(u_i$ simple and $u_n\neq 1$ unless $n=0)$ such that if also $u=u_1\cdots u_{k-1} v_1\cdots v_\ell$ $(v_i$ also simple, $\ell\geq 1)$ then $u_k=v_1w$ for some $w$.
\item The group $G$ is automatic.\footnote{See \cite{eps92} for a definition of automatic groups.}
\end{list}
\end{theorem}


\subsection{Complemented presentations} \label{ga95}

We turn to  Garside  groupoids.
Fix a finite set $G_0$ (which will be the object set of our categories). Let $A(x,y)$, $\barr{A}(x,y)$ ($x,y\in G_0$) be pairwise disjoint sets (which will provide the generating morphisms from $x$ to $y$). Let $a\mapsto a^{-1}$ denote a bijective map $A(x,y)\ra \barr{A}(y,x)$ (note the orders of $x$ and $y$); we shall denote its inverse by $a\mapsto a^{-1}$ too. The elements of $A(x,y)$ are called {\em letters}. Let $W$ be the category with object set $G_0$ and presented by the generators $A(x,y)\cup\barr{A}(x,y)$ (consisting of morphisms from $x$ to $y$) for all $x,y\in G_0$, and no relations. Let $W^+$ be the subcategory of $W$ with the same objects and generated by the letters. We call morphisms of $W$ {\em words\,} and morphisms of $W^+$ {\em positive words}. If $u$ is a word then we define $u^{-1}$ as the word obtained by inverting all letters occurring in $u$ and reversing their order. Note that $u^{-1}u\neq 1$ except if $u$ is a trivial morphism.

The (disjoint) union $\cup_y A(x,y)$ will be written $A(x,-)$. Let a positive word $f(a,b)$ be given for all $a,b\in A(x,-)$, $x\in G_0$. Suppose that, for all $a,b$, the positive words $a\,f(a,b)$ and $b\,f(b,a)$ are defined%
\footnote{The product $uv$ of two morphisms $u,v$ in a category is {\em defined\,} if the target object of $u$ is the source object of $v$.}
and have the same target object (they already have the same source object). The quotient of $W^+$ by the relations
\be \left\{a\,f(a,b)=b\,f(b,a)\ \Big|\ a,b\in A(x,-),\ x\in G_0\right\} \label{ga32} \ee
is written $G^+$. It is precisely the category presented by the same generators as $W^+$ and the relations (\ref{ga32}). The {\em groupoid\,} of the same presentation is written $G$. So there is a natural functor $G^+\ra G$. We write $u\equivp v$ ($u,v$ morphisms in $W^+$) if the images in $G^+$ of $u,v$ are equal.

We call \[ \big(G_0,\{A(x,y)\}_{xy},f\big) \] a {\em complemented presentation\,} or simply {\em complement}. In the rest of this section, we fix a complement and we use the above notation.

\subsection{The cube condition} \label{ga85}

For letters $a,b$ and words $u,v$ write
\[ u \, a^{-1} \, b \, v \ra u f(a,b) \, f(b,a)^{-1} v. \]
provided the left hand side is defined (in which case so is the right hand side). For words $u,v$ write $u\cra v$ if there are $n\geq 0$ and $u_i$ such that
\[ u=u_0\ra u_1\ra\cdots\ra u_n=v. \]

Suppose $u,v,u_0,v_0$ are positive words such that $u^{-1}v\cra v_0u_0^{-1}$. It can be shown that $v_0$ depends only on $u,v$; it will be written $u\un v$ and may or may not exist. Note $a\un b=f(a,b)$ for all $a,b\in A(x,{-})$.

We say that our complement satisfies the {\em cube condition\,} if
\be (a\un b)\un (a\un c)\equivp (b\un a)\un (b\un c) \label{ga33} \ee
for all $a,b,c\in A(x,-)$ (and all $x\in G_0$).

\subsection{Main theorem on Garside groupoids}

\begin{definition} \label{ga84} A category is called {\em atomic\,} if for each morphism $x$ there is a natural number $n$ such that $x$ cannot be written as a product of $n$ or more nontrivial morphisms.
\end{definition}

\begin{definition} By an {\em automorphism\,} of a complement $(G_0,A,f)$ (written on the right) we mean an automorphism $\phi$ of the associated category $W^+$ of positive words, such that $(A(x,y))\phi=A(x\phi,y\phi)$ for all $x,y\in G_0$ and $(f(a,b))\phi=f(a\phi,b\phi)$ whenever $a,b\in A(x,-)$. Note that the action of $\phi$ on $G_0$ may be nontrivial.
\end{definition}

\begin{definition} \label{ga71} A {\em Garside automorphism\,} is a pair \[ \big(\phi,\{\Delta_x\mid x\in G_0\}\big) \] where $\phi$ is an automorphism of the complement $(G_0,A,f)$ and $\Delta_x\in G^+(x,x\phi)$ such that the following hold.
\lista{{\rm(\GA\arabic{daana})}}{\setlength{\leftmargin}{4em}}
\item Whenever $a\in A(x,y)$ there exists $u\in G^+(y,x\phi)$ such that $au=\Delta_x$ (in $G^+$).
\item Whenever $a\in A(x,y)$ we have $a\Delta_y=\Delta_x(a\phi)$ (in $G^+$).
\end{list}
If $\Delta_x=uv$ (in $G^+$) then we call $u$ a {\em simple morphism}.
\end{definition}

\begin{definition} A {\em Garside tuple\,} consists of a complement $(G_0,A,f)$ and a Garside automorphism $(\phi,\{\Delta_x\}_x)$ such that $G^+$ is atomic, and the cube condition (\ref{ga33}) is satisfied.

In this case, we call $G$ a {\em Garside groupoid\,} and $G^+$ a {\em Garside category}.
\end{definition}

The following is the generalisation of theorem~\ref{ga101} to Garside groupoids, and the proof is easily adapted from the proof for the group case. We will use (\PC\ref{ga35}) and (\PC\ref{ga37}) later on.


\begin{theorem} \label{ga41} Let $(G_0,A,f,\phi,\{\Delta_x\}_x)$ be a Garside tuple. Then the following hold.
\lista{{\rm(\PC\arabic{daana})}}{\setlength{\leftmargin}{3em}}
\item \label{ga52} Let $[u]$ denote the image in $G^+$ of a morphism $u$ in $W^+$. If $u,v\in W^+(x,-)$ then $[u\un v]$ depends only on $[u]$ and $[v]$ and will therefore be written $[u]\un[v]$. 
\item \label{ga35} The category $G^+$ has finite limits and colimits. The limit $u\vee v$ of $\{u,v\}\subset G^+(x,{-})$ is $u(u\un v)=v(v\un u)$.
\item \label{ga37} The natural functor $G^+\ra G$ is injective.
\item \label{ga38} Every morphism in $G$ is a quotient of two morphisms in the image of the natural functor $G^+\ra G$.
\item \label{ga39} If $\Delta_x=uv$ in $G^+$ then $v$ is a simple morphism.
\item \label{ga40} {\rm (}Greedy form{\,\rm ).} Every morphism $u$ in $G^+$ can uniquely be written $u=u_1\cdots u_n$ $(u_i$ simple and $u_n\neq 1$ unless $n=0)$ such that whenever also $u=u_1\cdots u_{k-1} v_1\cdots v_\ell$ $(v_i$ also simple, $\ell\geq 1)$ then $u_k=v_1w$ for some $w$.
\item The groupoid $G$ is automatic.
\end{list}
\end{theorem}


%
      \section{The main result} \label{ga91}
%

\subsection{Introduction}

The elementary morphisms and elementary relations define a complemented groupoid presentation. In this section, we prove that it is part of a Garside tuple.

Let us recall what we already have. Let $G_0=M\backslash L$ be the object set of $G$ (definition~\ref{ga82}). Let $A(x,y)$ be the set of elementary morphisms from $x$ to $y$ ($x,y\in G_0$), see definition~\ref{ga3}. If $a_1,b_1\in A(x,y)$ then there is a unique elementary relation of the form $a_1\cdots a_k=b_1\cdots b_\ell$, see~\ref{ga83}. The elementary morphisms and relations yield a complemented presentation as defined in subsection~\ref{ga95}, which should not be confused with complementary generators and relations (subsection~\ref{ga72}). The complement $f$ is of course defined by $f(a_1,b_1):=a_2\cdots a_k$.

We merge the notations of the two previous sections. Note that each of the two sections gave their own definition of $G$, but that the two definitions are equivalent by proposition~\ref{ga28}.

In subsection~\ref{ga80} we prove atomicity (\ref{ga84}). In subsection~\ref{ga81} we prove the cube condition (\ref{ga85}). In the remaining subsections we construct a Garside automorphism (\ref{ga88}, \ref{ga57}) and prove that it is one.

\subsection{Atomicity} \label{ga80}

As every group, $\zz$ can be taken to be a groupoid. We will write $\zz_G$ for it to remind us that it is a groupoid.

It is well-known that there is a surjective homomorphism from the braid group to $\zz$. It follows that there is a surjective functor $w\col G\ra\zz_G$. In order to establish that it takes elementary morphisms to nonnegative integers, we construct $w$ in an independent way.

\begin{notation} \label{ga60} If $x$ is an admissible decomposition and $p$ a puncture we will write $R(x,p)$ for the region of $x$ containing~$p$. We will write $\arcs(x)$ for the set of arcs of~$x$.
\end{notation}

\begin{lemdef} \label{ga61} \rm (a). There exists a unique functor $w\col G\ra\zz_G$ with the following property. See figure~\ref{ga5}. Let $(x,y)$ be an elementary pair. Write 
\begin{align*} \arcs(x)\smallsetminus\arcs(y) &=\{A\}, \\  \arcs(y)\smallsetminus\arcs(x) &=\{B\}, \\ A&=R(x,p)\cap R(x,q) \end{align*} so that also $B=R(y,p)\cap R(y,q)$. Then $w$ takes the elementary morphism $M(x,y)$ to the number of boundary edges of $D$ which are separated from $p$ by both or neither of $A,B$. As an example, we have indicated the boundary edges with this property by $w$ in figure~\ref{ga5}; there are $6$ of them so $w$ takes the elementary morphism of the figure to~$6$.

(b). The functor $w\col G\ra\zz_G$ takes elementary morphisms to nonnegative integers. 

(c). The category $G^+$ is atomic.
\end{lemdef}

\begin{proof} (a). One proves existence of $w$ by checking that $w$ takes every elementary relation to a relation in $\zz_G$; this can quickly be done using the diagrams for the elementary relations in figure~\ref{ga26}. 

(b). This is obvious.

(c). Let $G^+_0$ denote the subcategory of $G^+$ generated by those elementary morphisms $x$ satisfying $w(x)=0$. By (b) it is enough to prove that $G^+_0$ is atomic.

Consider first the case where all labels are~$2$. Then $G^+_0=G^+$. Also, $G^+$ is isomorphic to the positive braid monoid $B_n^+\subset B_n$ generated by $\sigma_1$, \ldots, $\sigma_{n-1}$, which is well-known to be atomic; this follows for example from the existence of a homomorphism $B_n\ra\zz$ which takes each $\sigma_i$ to~$1$.

Finally we consider the general case. Notice that every morphism of $G^+_0$ is an endomorphism. It is therefore enough to prove the atomicity of $\End(X)$, the monoid of endomorphisms in $G^+_0$ of any object~$X$. Now $\End(X)$ is isomorphic to a direct product of positive braid monoids which is atomic by the foregoing.\end{proof}

\subsection{Proof of the cube condition} \label{ga81}

\begin{convention} \label{ga73} Recall our convention~\ref{ga30} on the meaning of apparent edges in objects. In addition, from now on, thin apparent edges will stand for precisely one edge.
\end{convention}

\begin{lemma} \label{ga70} The cube condition {\rm (\ref{ga85})} is satisfied.
\end{lemma}

\begin{proof} Let us start with an example. See figure~\ref{ga42}(a). We will prove 
\be (a\un b)\un (a\un c)\equivp (b\un a)\un (b\un c), \label{ga64} \ee
where $a,b,c$ are defined by the figure, and which is an instance of the cube condition (\ref{ga33}). We have
\begin{align*}
&a\un b=d,\quad a\un c=k\ell,\quad b\un a=ef,\quad  b\un c= g, \\
&(a\un b)^{-1}(a\un c)=\underline{d^{-1}\ k}\,\ell \cra n\,\underline{p^{-1}\,\ell} \cra n\,q\,r^{-1}, \\
&(b\un a)^{-1}(b\un c)= f^{-1}\,\underline{e^{-1}\,g} \cra \underline{f^{-1}\,s}\,t^{-1} \cra n\,q\,u^{-1}\,t^{-1}
\end{align*}
so
\begin{align*} (b\un a)\un(b\un c)=nq=(a\un b)\un(a\un c). \end{align*}
It is not hard to do this calculation without writing much down, if you have figure~\ref{ga42}(a) in front of you. As the same diagram helps you verify two more equations (namely (\ref{ga64}) with $a,b,c$ permuted), drawing these diagrams is a useful step towards verifying the cube condition.

I haven't explained yet how the diagrams come about, even though we have rigorously used one already.

\newcommand{\eightgon}{%
\pspicture(0,0)(0,0) \SpecialCoor \psset{unit=5mm} \degrees[80]
\pnode(1; 0){q0} \pnode(1;10){q1}
\pnode(1;20){q2} \pnode(1;30){q3}
\pnode(1;40){q4} \pnode(1;50){q5}
\pnode(1;60){q6} \pnode(1;70){q7}
\psarc(0,0){1}{(q0)}{(q1)} \psarc(0,0){1}{(q2)}{(q3)}
\psarc(0,0){1}{(q4)}{(q5)} \psarc(0,0){1}{(q6)}{(q7)}
\psset{linestyle=solid, linewidth=1.8pt, dash=2.4pt 1.5pt}
\psarc(0,0){1}{(q1)}{(q2)} \psarc(0,0){1}{(q3)}{(q4)}
\psarc(0,0){1}{(q5)}{(q6)} \psarc(0,0){1}{(q7)}{(q0)}
\endpspicture}
\newcommand{\eightdots}{{\SpecialCoor%
\psdots(q0)(q1)(q2)(q3)(q4)(q5)(q6)(q7)}}
\newcommand{\sixgon}{\pspicture(0,0)(0,0) \degrees[12] \psset{unit=3.5mm} \pnode(1;3){q0}  \pnode(1;5){q1}  \pnode(1;7){q2}  \pnode(1;9){q3}  \pnode(1;11){q4}  \pnode(1;13){q5} \pscircle(0,0){1} \endpspicture}
\newcommand{\sixdots}{\psdots(q0)(q1)(q2)(q3)(q4)(q5)}
\newcommand{\dron}[3]{\psframebox[linestyle=none, framesep=3mm]{\pspicture(-24mm, -29mm)(24mm, 29mm)
\SpecialCoor \psset{xunit=.75mm, yunit=.75mm, dimen=middle, dotstyle=o, dotsize=3pt 0, linewidth=.8pt}
\rput(-20,32){\text{\footnotesize (#3)}}
\pnode(-20,0)  {n1}
\pnode(20,0)   {n5}
\pnode(0,10)  {n2}
\pnode(0,-10)   {n3}
\pnode(-20,20){n7}
\pnode(20,20) {n8}
\pnode(20,-20)  {n9}
\pnode(-20,-20) {n14}
\pnode(0,30)  {n4}
\pnode(0,-30)   {n6}
\pnode(-32,38){n10}
\pnode(32,38) {n11}
\pnode(32,-38)  {n12}
\pnode(-32,-38) {n13}
{\psset{fillstyle=solid, fillcolor=lightgray, linestyle=none} \pscircle(#1){3.5mm} \pscircle(#2){3.5mm}}
\rput(n1){\sixgon} \psline(q5)(q1)(q4) \psline(q1)(q3) \sixdots
\rput(n2){\sixgon} \psline(q2)(q4)(q1)(q5) \sixdots
\rput(n3){\sixgon}  \psline(q4)(q2)(q5)(q1) \sixdots
\rput(n4){\sixgon}  \psline(q0)(q4)(q2)  \psline(q1)(q4) \sixdots
\rput(n5){\sixgon}  \psline(q4)(q2)(q0) \psline(q5)(q2) \sixdots
\rput(n6){\sixgon}  \psline(q2)(q5)  \psline(q3)(q5)(q1) \sixdots
\rput(n7){\sixgon}  \psline(q3)(q1)(q4)(q0) \sixdots
\rput(n8){\sixgon}  \pspolygon(q4)(q2)(q0) \sixdots
\rput(n9){\sixgon} \psline(q3)(q5)(q2)(q0)  \sixdots
\rput(n10){\sixgon}  \psline(q4)(q0)(q3)(q1) \sixdots
\rput(n11){\sixgon}  \psline(q3)(q0) \psline(q4)(q0)(q2) \sixdots
\rput(n12){\sixgon} \psline(q5)(q3)(q0)(q2)  \sixdots
\rput(n13){\sixgon} \psline(q0)(q3) \psline(q5)(q3)(q1)  \sixdots
\rput(n14){\sixgon} \pspolygon(q5)(q3)(q1)  \sixdots
%
\endpspicture}}

\begin{figure}
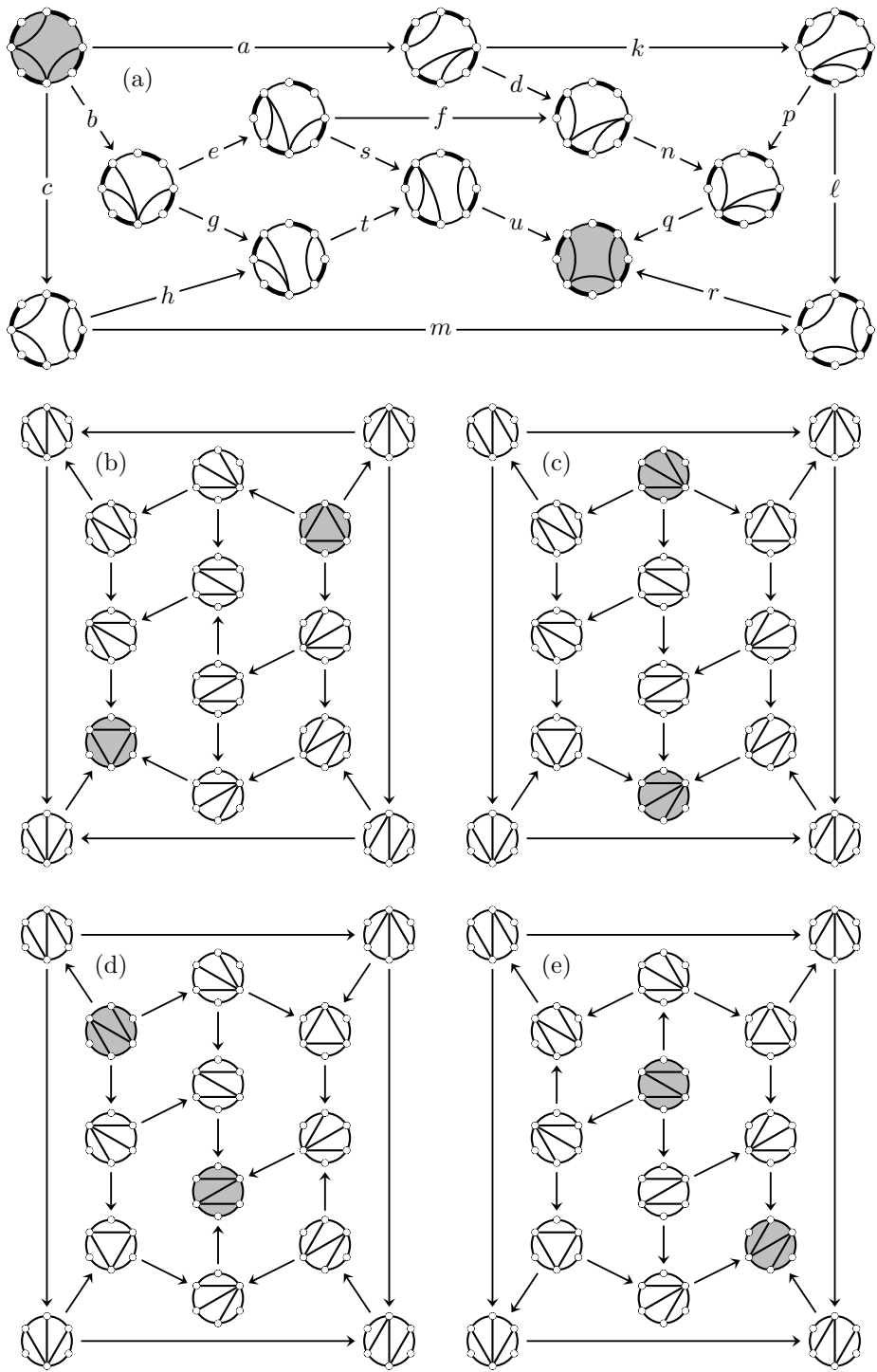
 \footnotesize 
\centering
{\psset{arrowsize=2pt 4, arrowlength=.7, linewidth=.7pt, framesep=.4ex}
\psframebox[linestyle=none, framesep=5mm]{%
\pspicture(-30mm,1mm)(30mm,40mm) $
\psset{xunit=.85mm, yunit=.66mm, dotstyle=o, dotsize=3.5pt, linewidth=.8pt}
\rput(-50,53){\text{\footnotesize (a)}}
\pnode(-65,  60){n1}  \pnode(   0,  60){n2}
\pnode(-65,   0){n3}  \pnode( -50,  30){n4}
\pnode( 65,   0){n5}  \pnode(  65,  60){n6}
\pnode(-25,  45){n7}  \pnode(  25,  45){n8}
\pnode(-25,  15){n9}  \pnode(  50,  30){n10}
\pnode(  0,  30){n11} \pnode(  25,  15){n12}
\psset{labelsep=0, arcangle=40}
\newpsobject{myarc}{ncarc}{arcangle=18}
{\SpecialCoor\psset{fillstyle=solid, fillcolor=lightgray, linestyle=none} \pscircle(n1){5mm} \pscircle(n12){5mm}}
\nput{0}{n1}{\eightgon} \ncarc{q2}{q4} \ncarc{q4}{q6} \ncarc{q6}{q0} \eightdots
\nput{0}{n2}{\eightgon} \ncarc{q2}{q4} \ncarc{q6}{q0} \myarc{q5}{q0} \eightdots
\nput{0}{n3}{\eightgon} \ncarc{q2}{q4} \ncarc{q4}{q6} \ncarc{q7}{q1} \eightdots
\nput{0}{n4}{\eightgon} \myarc{q3}{q6} \ncarc{q4}{q6} \ncarc{q6}{q0} \eightdots
\nput{0}{n5}{\eightgon} \ncarc{q2}{q4} \ncarc{q5}{q7} \ncarc{q7}{q1} \eightdots
\nput{0}{n6}{\eightgon} \ncarc{q2}{q4} \ncarc{q5}{q7} \myarc{q5}{q0} \eightdots
\nput{0}{n7}{\eightgon} \ncarc{q3}{q5} \myarc{q3}{q6} \ncarc{q6}{q0} \eightdots
\nput{0}{n8}{\eightgon} \ncarc{q3}{q5} \myarc{q5}{q0} \ncarc{q6}{q0} \eightdots
\nput{0}{n9}{\eightgon} \myarc{q3}{q6} \ncarc{q4}{q6} \ncarc{q7}{q1} \eightdots
\nput{0}{n10}{\eightgon} \ncarc{q3}{q5} \ncarc{q5}{q7} \myarc{q5}{q0} \eightdots
\nput{0}{n11}{\eightgon} \ncarc{q3}{q5} \myarc{q3}{q6} \ncarc{q7}{q1} \eightdots
\nput{0}{n12}{\eightgon} \ncarc{q3}{q5} \ncarc{q5}{q7} \ncarc{q7}{q1} \eightdots
\psset{arrowsize=2pt 4, arrowlength=.7, linewidth=.7pt, arrows=->, nodesep=6.4mm}
\ncline{n1}{n2}   \ncput*{a}
\ncline{n1}{n4}   \ncput*{b}
\ncline{n2}{n6}   \ncput*{k}
\ncline{n2}{n8}   \ncput*{d}
\ncline{n3}{n5}   \ncput*{m}
\ncline{n3}{n9}   \ncput*{h}
\ncline{n4}{n7}   \ncput*{e}
\ncline{n4}{n9}   \ncput*{g}
\ncline{n5}{n12}  \ncput*{r}
\ncline{n7}{n8}   \ncput*{f}
\ncline{n7}{n11}  \ncput*{s}
\ncline{n8}{n10}  \ncput*{n}
\ncline{n9}{n11}  \ncput*{t}
\ncline{n10}{n12} \ncput*{q}
\ncline{n11}{n12} \ncput*{u}
\psset{framesep=.5ex}
\ncline{n1}{n3}   \ncput*{c}
\ncline{n6}{n5}   \ncput*{\ell}
\ncline{n6}{n10}  \ncput*{p}
$\endpspicture} \\[4ex]
\psset{nodesep=4.8mm}
$\begin{array}{@{}c@{\hspace{8mm}}c@{}}
\dron{n8}{n14}{b}
\ncline{<-}{n1}{n2} \ncline{<-}{n1}{n7} \ncline{->}{n1}{n14} \ncline{<-}{n2}{n3} \ncline{<-}{n2}{n4} \ncline{<-}{n3}{n5} \ncline{->}{n3}{n6} \ncline{->}{n4}{n7} \ncline{<-}{n4}{n8} \ncline{<-}{n5}{n8} \ncline{->}{n5}{n9} \ncline{<-}{n6}{n9} \ncline{->}{n6}{n14} \ncline{->}{n7}{n10} \ncline{->}{n8}{n11} \ncline{<-}{n9}{n12} \ncline{<-}{n10}{n11} \ncline{->}{n10}{n13} \ncline{->}{n11}{n12} \ncline{->}{n12}{n13} \ncline{->}{n13}{n14}
&
\dron{n4}{n6}{c}
\ncline{<-}{n1}{n2} \ncline{<-}{n1}{n7} \ncline{->}{n1}{n14} \ncline{->}{n2}{n3} \ncline{<-}{n2}{n4} \ncline{<-}{n3}{n5} \ncline{->}{n3}{n6} \ncline{->}{n4}{n7} \ncline{->}{n4}{n8} \ncline{<-}{n5}{n8} \ncline{->}{n5}{n9} \ncline{<-}{n6}{n9} \ncline{<-}{n6}{n14} \ncline{->}{n7}{n10} \ncline{->}{n8}{n11} \ncline{<-}{n9}{n12} \ncline{->}{n10}{n11} \ncline{->}{n10}{n13} \ncline{->}{n11}{n12} \ncline{<-}{n12}{n13} \ncline{->}{n13}{n14} \\[8mm]
\dron{n7}{n3}{d}
\ncline{->}{n1}{n2} \ncline{<-}{n1}{n7} \ncline{->}{n1}{n14} \ncline{->}{n2}{n3} \ncline{<-}{n2}{n4} \ncline{<-}{n3}{n5} \ncline{<-}{n3}{n6} \ncline{<-}{n4}{n7} \ncline{->}{n4}{n8} \ncline{<-}{n5}{n8} \ncline{<-}{n5}{n9} \ncline{<-}{n6}{n9} \ncline{<-}{n6}{n14} \ncline{->}{n7}{n10} \ncline{<-}{n8}{n11} \ncline{<-}{n9}{n12} \ncline{->}{n10}{n11} \ncline{->}{n10}{n13} \ncline{->}{n11}{n12} \ncline{<-}{n12}{n13} \ncline{->}{n13}{n14}
&
\dron{n2}{n9}{e}
\ncline{<-}{n1}{n2} \ncline{->}{n1}{n7} \ncline{->}{n1}{n14} \ncline{->}{n2}{n3} \ncline{->}{n2}{n4} \ncline{->}{n3}{n5} \ncline{->}{n3}{n6} \ncline{->}{n4}{n7} \ncline{->}{n4}{n8} \ncline{<-}{n5}{n8} \ncline{->}{n5}{n9} \ncline{->}{n6}{n9} \ncline{<-}{n6}{n14} \ncline{->}{n7}{n10} \ncline{->}{n8}{n11} \ncline{<-}{n9}{n12} \ncline{->}{n10}{n11} \ncline{->}{n10}{n13} \ncline{->}{n11}{n12} \ncline{<-}{n12}{n13} \ncline{<-}{n13}{n14}
\end{array}$}  
\caption{The irreducible rank 3 characteristic graphs if all labels are~$\geq 3$. \label{ga42}}
\end{figure}

\begin{definition} A commutative diagram over $G^+$ all of whose arrows are elementary morphisms is {\em closed\,} if the following holds for every  elementary relation $a_1\cdots a_k=b_1\cdots b_\ell$. If the diagram contains $a_1$ and $b_1$ then it contains the whole elementary relation.

Let $I\subset A(x,-)$ be any subset. The {\em characteristic graph spanned by $I$} is the (unique) closed diagram $\CG(I)$ containing $I$ and which is universal for this property. The {\em rank\,} of the characteristic graph is $|I|$ and its {\em initial object\,} is~$x$.
\end{definition}

Verifying the cube condition is done in two steps: one, drawing all rank~3 characteristic graphs; two, verifying the cube condition for each of them as we did in our example of figure~\ref{ga42}(a).

\begin{remark} \label{ga96} Assume for simplicity that every characteristic graph is finite. It can be shown that then the cube condition holds if and only if every characteristic graph has limits. We shall not use this.
\end{remark}

\begin{definition} Let $A$ be an arc of an object $x\in G_0$. The arc $A$ cuts the disk $D$ into two pieces, say, $D_1$ and $D_2$. Let $I\subset A(x,-)$ be such that $A$ is not the rotating arc in any element of~$I$. Then $I$ is a disjoint union $I=I_1\amalg I_2$ where $I_i$ takes place in~$D_i$. 

The characteristic graph $\CG(I)$ is said to be {\em reducible\,} if (for some choice of $A$) the above properties hold and $I_1,I_2$ are non-empty. We say it is {\em complete\,} if $I=A(x,-)$.
\end{definition}

It is easy to prove the cube condition in reducible rank~3 characteristic graphs. We turn to the irreducible ones. 

By the {\em shape\,} of a characteristic graph we mean the underlying directed graph up to isomorphism. Figure~\ref{ga42}(a) gives one shape in detail. Just as the elementary relations (ER2--4), it involves fat apparent edges which stand for unspecified nonnegative numbers of edges, the shape of the characteristic graph being independent of the exact numbers.

Checking the cube condition on a characteristic graph involves only the shape of that graph. Every irreducible characteristic graph has the shape of a complete one. So we need only look at complete rank~3 characteristic graphs from now on. They are completely determined by their initial objects, which are elements of $G_0$ of $3$~arcs.

Figure~\ref{ga63} gives the initial objects of all (complete rank~$3$) characteristic graphs partitioned in eight {\em cases\,} (a)--(h) with the property that the characteristic graphs are of the same shape if (and only if) they are in the same case, as one finds when one draws these graphs explicitly, and which the reader should do. The figure also gives the number of objects in the characteristic graphs. Note that the boundary edges of the initial object can affect the shape of the graph; for example, by contracting two appropriate apparent boundary edges of (a) one obtains (b) which represents a different shape.

The characteristic graphs of cases (a)--(e) are drawn in figure~\ref{ga42}; these are precisely the cases which occur if all labels are $\geq 3$. For convenience, the initial and terminal objects are gray.
Case (a) is the only rank~$3$ shape that does not occur if all labels are $\leq 3$. For simplicity the pictures of the characteristic graphs (b)--(e) assume that all labels are~3, that is, they avoid the fat apparent edges. Using figure~\ref{ga42} and following our example of case~(a) it is then straightforward and not-so-tedious to verify the cube condition for the case where all labels are at least~$3$.

As to cases (f), (g), (h), we leave it to the reader to draw their characteristic graphs and to verify the cube condition.

This finishes our proof of the cube condition.\end{proof}

\newcommand{\dashc}{\psset{linewidth=1.8pt}}

\newcommand{\fourdotsc}{\psdots(q0)(q1)(q2)(q3)}
\newcommand{\sixdotsc}{\psdots(q0)(q1)(q2)(q3)(q4)(q5)}
\newcommand{\eightdotsc}{\psdots(q0)(q1)(q2)(q3)(q4)(q5)(q6)(q7)}

\newcommand{\fourgonc}{%
{\pspicture[.5](-6mm,-6mm)(6mm,6mm) \psset{unit=.06mm, dimen=middle} \degrees[4]
\pnode(100;1){q0} \pnode(100;2){q1} \pnode(100;3){q2} \pnode(100;4){q3}
\psarc(0,0){100}{(q0)}{(q1)} \psarc(0,0){100}{(q2)}{(q3)} 
\dashc \psarc(0,0){100}{(q1)}{(q2)} \psarc(0,0){100}{(q3)}{(q0)} 
\endpspicture}}

\newcommand{\sixgonc}{%
{\pspicture[.5](-6mm,-6mm)(6mm,6mm) \psset{unit=.06mm, dimen=middle} \degrees[12]
\pnode(100;3){q0} \pnode(100;5){q1} \pnode(100;7){q2} \pnode(100;9){q3} \pnode(100;11){q4} \pnode(100;13){q5}
\psarc(0,0){100}{(q0)}{(q1)} \psarc(0,0){100}{(q2)}{(q3)} \psarc(0,0){100}{(q4)}{(q5)} 
\dashc \psarc(0,0){100}{(q1)}{(q2)} \psarc(0,0){100}{(q3)}{(q4)} \psarc(0,0){100}{(q5)}{(q0)}
\endpspicture}}

\newcommand{\eightgonc}{%
{\pspicture[.5](-6mm,-6mm)(6mm,6mm) \psset{unit=.06mm, dimen=middle} \degrees[8]
\pnode(100;2){q0} \pnode(100;3){q1} \pnode(100;4){q2} \pnode(100;5){q3} \pnode(100;6){q4} \pnode(100;7){q5} \pnode(100;8){q6} \pnode(100;1){q7}
\psarc(0,0){100}{(q0)}{(q1)} \psarc(0,0){100}{(q2)}{(q3)} \psarc(0,0){100}{(q4)}{(q5)} \psarc(0,0){100}{(q6)}{(q7)} 
\dashc \psarc(0,0){100}{(q1)}{(q2)} \psarc(0,0){100}{(q3)}{(q4)} \psarc(0,0){100}{(q5)}{(q6)} \psarc(0,0){100}{(q7)}{(q0)}
\endpspicture}}


\begin{figure}
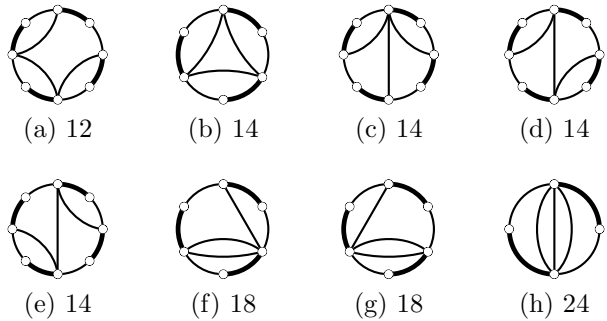
 
\centering
\SpecialCoor
\psset{dotstyle=o, dotsize=3.6pt, linewidth=.8pt, arcangle=40}
\footnotesize
\psmatrix[colsep=10mm, rowsep=2mm]
\eightgonc \ncarc{q0}{q2} \ncarc{q2}{q4} \ncarc{q4}{q6} \eightdotsc &
\sixgonc {\psset{arcangle=20} \ncarc{q0}{q2} \ncarc{q2}{q4} \ncarc{q4}{q0}} \sixdotsc &
\eightgonc \ncarc{q0}{q2} \ncline{q0}{q4} \ncarc{q6}{q0} \eightdotsc &
\eightgonc \ncarc{q0}{q2} \ncline{q0}{q4} \ncarc{q4}{q6} \eightdotsc \\
\text{(a) 12} & \text{(b) 14} & \text{(c) 14} & \text{(d) 14} \\[3mm]
\eightgonc \ncarc{q2}{q4} \ncline{q0}{q4} \ncarc{q6}{q0} \eightdotsc &
\sixgonc \ncarc[arcangle=40]{q2}{q4} \ncarc[arcangle=15]{q4}{q2} \ncline{q0}{q4} \sixdotsc &
\sixgonc \ncarc[arcangle=40]{q2}{q4} \ncarc[arcangle=15]{q4}{q2} \ncline{q0}{q2} \sixdotsc &
\fourgonc \psset{arcangle=50} \ncarc{q2}{q0} \ncarc{q0}{q2} \ncline{q0}{q2} \fourdotsc \\
\text{(e) 14} & \text{(f) 18} & \text{(g) 18} & \text{(h) 24} 
\endpsmatrix 
\caption{Initial objects of characteristic graphs classified by shape. \label{ga63}}
\end{figure}

\subsection{Shifted decompositions}

The remainder of this section, which involves several subsections, is devoted to a proof of the existence of a Garside automorphism. 

Let me begin by explaining what is difficult about it. Recall that we are aiming to construct $\Delta_x\in G^+(x,x\phi)$ satisfying two identities in the category $G^+$ stated in (\GA 1) and (\GA 2) in definition~\ref{ga71}. It is easy to establish these identities in $G$ rather than $G^+$; as the example $\<a,b\mid ab=bb\>$ shows, this is not enough though. We need to prove the identities to hold in $G^+$, and this is harder.

Let $D_0=\{z\in\cc:|z|\leq 1\}$. Fix a finite non-empty set $P$ of interior points of $D_0$ and a labelling $\ell\col P\ra\zz_{\geq 2}$.

Define $m$ by
\[ m-2=\sum_{p\in P}\big( \ell(p)-2\big), \]
compare (\ref{ga1}).
Define $S\col \partial D_0\ra\partial D_0$ by $S(z)=z\,\exp(2\pi i/m)$. By an $S$-orbit we mean a set of the form $\{S^k(q)\mid k\in\zz\}$ for some $q\in \partial D_0$.

\begin{definition} \label{ga43} A {\em shifted decomposition\,} consists of a closed disk $D\subset D_0$ such that $D\cap\partial D_0$ is an $S$-orbit  and such that the interior of $D$ contains $P$, together with an admissible decomposition for $(D,P,D\cap\partial D_0,\ell)$.
\end{definition}

We call $D\cap\partial D_0$ the vertex set of the shifted decomposition. We will soon be looking at pairs of shifted decompositions of distinct vertex sets. A {\em single\,} shifted decomposition is essentially the same as an admissible decomposition. Most of our conventions about admissible decompositions easily translate to shifted decompositions but we mention the following. 

An arc of a shifted decomposition involving a disk $D\subset D_0$ is by definition an edge not contained in $\partial D$ (it never is in $\partial D_0$).

Two shifted decompositions of the same vertex set are said to be {\em isotopic\,} if they are isotopic relative $\partial D_0$ (pointwise) and $P$. In other words, two shifted decompositions of the same vertex set may be isotopic but not involve the same $D$.

Pictures of shifted decomposition can be simplified by not drawing the boundary of $D$. One only draws $D\cap\partial D_0$ and the arcs, that is, interior edges. See figure~\ref{ga49} (d), (e) for an example.

If $Q$ is an $S$-orbit we will write $L(Q)$ for the set of isotopy classes of shifted decompositions with vertex set~$Q$. There is a natural bijection between $L$ and $L(S^\zz (1))$ where $S^\zz (1)=\{ S^k(1)\mid k\in\zz\}$ is the $S$-orbit of $1$, which is the set of complex $m$-th roots of unity.

Every $S$-orbit $Q$ yields $G(Q)$ and $G^+(Q)$ which are versions of $G$ and $G^+$ in an obvious way.

\renewcommand{\erfivespec}{%
\pspicture(0,0)(0,0) \SpecialCoor  \psset{unit=.085mm} \degrees[20] \pnode(100; 3){q1} \pnode(100; 7){q2} \pnode(100;11){q3} \pnode(100;15){q4} \pnode(100; 19){q0} \pscircle[dimen=middle](0,0){100} 
\psdots(30;15)(65;1)(65;9)
\psset{linestyle=dashed, dash=3pt 2pt, showpoints=true, radius=1.7pt}
\Cnode*(100;17){p1} \Cnode*(100;5){p2} \Cnode*(100;13){p3} \ncline{p1}{p2} \ncline{p2}{p3}
\psdots(100;1) \psdots(100;9)
\endpspicture}

\renewcommand{\erfivedots}{{\psdots[dotstyle=o](q0)(q1)(q2)(q3)(q4)}}

\renewcommand{\simplefive}{
\SpecialCoor \linea \psset{dotsize=3.4pt, dotstyle=*}
\psframebox[framesep=8.5mm, linestyle=none]{\pspicture[.5](0mm,0mm)(56mm, 22mm)  \psset{unit=2.8mm, arcangle=25, linewidth=.7pt}
\pnode( 0,8){n1}  \pnode(10,8){n2}
\pnode(20,8){n3}  \pnode( 5,0){n4}
\pnode(15,0){n5}
\rput(n1){\erfivespec} \ncline{q3}{q1} \ncline{q1}{q4}  \erfivedots
\rput(n2){\erfivespec}\ncline{q1}{q4} \ncline{q4}{q2} \erfivedots
\rput(n3){\erfivespec}\ncline{q4}{q2}  \ncline{q2}{q0} \erfivedots
\rput(n4){\erfivespec}\ncarc{q3}{q0} \ncline{q3}{q1}  \erfivedots
\rput(n5){\erfivespec}\ncarc{q3}{q0} \ncline{q2}{q0} \erfivedots
\psset{arrows=->, nodesep=9.5mm}
\ncline{n1}{n2} \ncline{n2}{n3} \ncline{n5}{n3} \ncline{n1}{n4} \ncline{n4}{n5}
\endpspicture}}

\subsection{Compatibility, part 1}

Let $x,y$ be two shifted decompositions of distinct vertex sets. After an isotopy on one of them we may assume that any edge of $x$ intersects any edge of $y$ transversally in a minimal number of points. The pair $(x,y)$ is then said to be {\em tight}. The simultaneous isotopy class $[(x,y)]$ of a tight pair depends only on the individual isotopy classes $[x]$ and $[y]$.

\begin{definition} \label{ga44} (See figure~\ref{ga49}). Let $x,y$ be two shifted decompositions of distinct vertex sets in tight position. We say that $[x],[y]$ are {\em compatible\,} if the following holds for all punctures~$p$. Loosely speaking, some rotation with centre $p$ takes $R(x,p)$ to $R(y,p)$. More precisely, the vertex sets of $R(x,p)$ and $R(y,p)$ alternate along the boundary; and there exists a homeomorphism $h\col D_0\ra D_0$ such that $hR(x,p)$ and $hR(y,p)$ are convex. (If $\ell(p)>2$ then an equivalent condition is that $hR(x,p)$ is the convex hull of $hR(x,p)\cap\partial D_0$ and the same for $y$).
\newboolean{BoundaryEdges}
\newcommand{\compat}{%
\footnotesize 
\degrees[240] \psset{unit=.20mm, arcangle=43, linewidth=.6pt}
\psframebox[linestyle=none, framesep=0]{%
\pspicture(-100,-100)(100,100)$ 
\pscircle(0,0){100}
\ifthenelse{\boolean{BoundaryEdges}}{
\multido{\ia=0+20, \ib=20+20}{24}{\pcarc(100;\ia)(100;\ib)}}{}
\pcarc[arcangle=30](100;20)(100;80)
\pcarc[arcangle=30](100;80)(100;140)
\pcarc[arcangle=40](100;140)(100;180)
\pnode(100;180){q1} \pnode(100;20){q2} \pnode(100;200){q3} \pnode(100;80){q4}
\pcarc[arcangle=40](100;220)(100;20)
\pcarc[arcangle=20](q1)(q2)
\psset{linestyle=dashed, dash=1.5pt 1.2pt}
\ifthenelse{\boolean{BoundaryEdges}}{%
\multido{\ia=10+20, \ib=30+20}{24}{\pcarc(100;\ia)(100;\ib)}}{}
\pcarc[arcangle=30](100;30)(100;90)
\pcarc[arcangle=30](100;90)(100;150)
\pcarc[arcangle=30](100;170)(100;230)
\pcarc[arcangle=20](100;150)(100;230)
\pcarc[arcangle=12](100;150)(100;10)
\multido{\ia=0+20}{12}{\psdot(100;\ia)}
{\psset{dotstyle=o}
\multido{\ia=10+20}{12}{\psdot(100;\ia)}}
\psset{radius=1.7pt, labelsep=3pt}
\Cnode*(-10,20){a} \nput{-60}{a}{4}
\Cnode*(5,65){a} \nput{0}{a}{4}
\Cnode*(-65,15){a} \nput{-60}{a}{4}
\Cnode*(75;160){a} \nput{-20}{a}{3}
\Cnode*(75,3){a} \nput{-60}{a}{3}
\Cnode*(35,-57){a} \nput{20}{a}{4}
$\endpspicture}}
\newcommand{\compattwo}{%
\psframebox[linestyle=none, framesep=0]{%
\psset{unit=.10mm, linestyle=solid, radius=1.7pt, labelsep=3pt}
\pspicture[.85](-100,-100)(100,100)
\pscircle[linewidth=.6pt](0,0){100}
\degrees[4]
\pnode(100;0){q0} \pnode(100;1.5){q1} \pnode(100;2){q2} \pnode(100;3.5){q3}
\psline[linestyle=dashed, dash=3pt 2pt](q1)(q3)
\pnode(-20,60){p1} \pnode(20,60){p2} \pnode(-20,-60){p3} \pnode(20,-60){p4} 
\psdots[dotstyle=o](q1)(q3) \psdots(q0)(q2)
\psset{radius=1.7pt}
\Cnode*(0,39){a} \Cnode*(0,-39){a}
\endpspicture}}%
\newcommand{\compatc}{%
\begin{array}{@{}c@{}} 
\simplefive \\ \text{(c).\ \rule{0mm}{4ex}} \parbox[t]{67mm}{A shifted decomposition $y\in L(Q_2)$ (dashed) and all elements of $L(Q_1)$ compatible with it (solid)}  
\end{array}}
\newcommand{\compatab}{%
\begin{array}{@{}c@{}}
\compattwo \psline(q0)(q2) \\ \text{(a). Compatible\rule{0mm}{2.5ex}} \\ \\
\compattwo \pscurve(q0)(p3)(p2)(q2) \\ \text{(b). Not compatible\rule{0mm}{2.5ex}}
\end{array}}
\newcommand{\compatd}{%
\begin{array}{@{}c@{\hspace{10mm}}c@{}}
\setboolean{BoundaryEdges}{false} \compat 
& \setboolean{BoundaryEdges}{true} \compat \\[2mm]
\text{\normalsize (d). Compatible} & \text{\normalsize (e).\ \parbox[t]{44mm}{\normalsize The same with the boundary edges drawn}}
\end{array}}
\begin{figure}
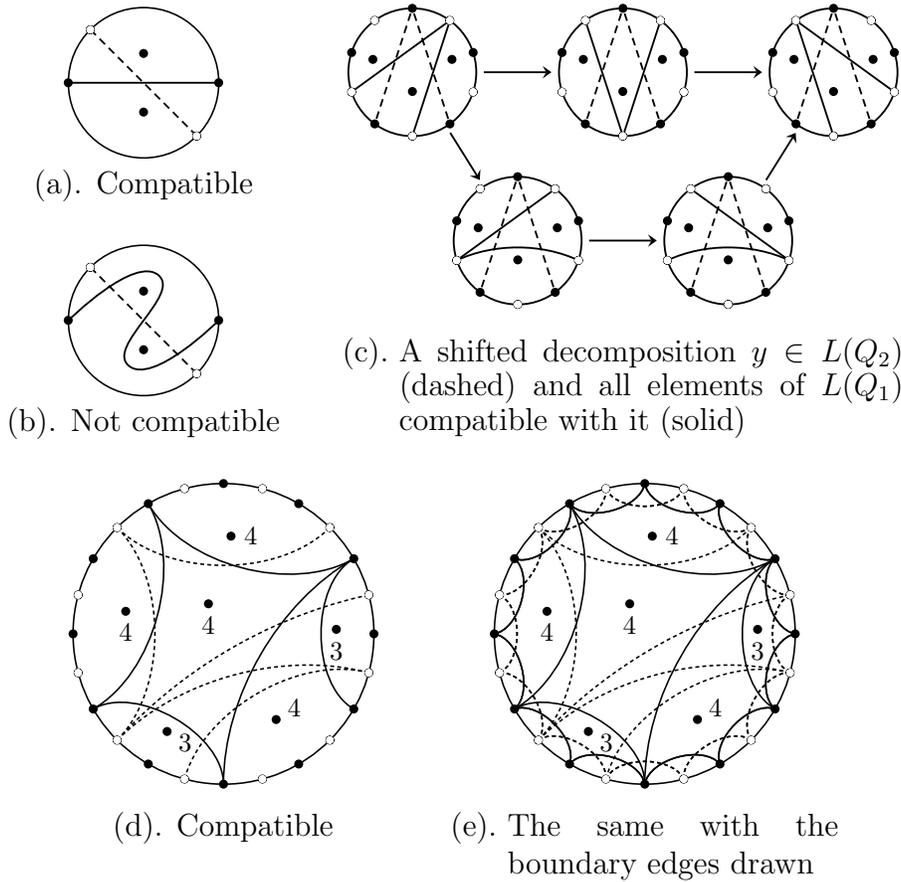
 \SpecialCoor \psset{dotsize=3.4pt 0, linewidth=.7pt, dimen=middle} 
\begin{gather*}
\compatab \qquad \compatc \\[2ex] \compatd
\end{gather*}
\caption{Compatibility. See definition~\ref{ga44}. Elements in $L(Q_1)$ are solid and those in $L(Q_2)$ dashed. \label{ga49}}
\end{figure}%
\end{definition}

\begin{definition}
From now on we fix two distinct $S$-orbits $Q_1$, $Q_2$. For every $y\in L(Q_2)$ we will write $\Omega(y)$ for the set of those $x\in L(Q_1)$ compatible with~$y$. Elements of $L(Q_1)$ are drawn solid and those of $L(Q_2)$ dashed.
\end{definition}

\begin{example} If the number of arcs is $n-1=2$ then the set of objects in an elementary relation (or rather their lifts to $L(Q_1)$) is of the form $\Omega(y)$. See figure~\ref{ga49}(c) where $y$ is indicated by dashed arcs. Similarly, each of the diagrams in figure~\ref{ga42} (or again rather their lifts to $L(Q_1)$) is of the form $\Omega(y)$, if $n-1=3$.
\end{example}

The following two lemmas are easy and left to the reader to prove.

\begin{lemma} \label{ga45} Let $([x],[y])\in L(Q_1)\times L(Q_2)$ be compatible and suppose that $x,y$ are tight. \begin{enumerate}
\item Any edge of $x$ meets any edge of $y$ in at most one point.
\item The set $\Omega([y])$ is finite.\qed
\end{enumerate}
\end{lemma}

From now on we will abuse language as usual and confuse shifted decompositions with their isotopy classes.

\begin{lemdef} \rm \label{ga59} Let $(x_2,y)\in L(Q_1)\times L(Q_2)$ be compatible. Let $A$ be an arc of $x_2$, and let $(x_1,x_2)$, $(x_2,x_3)$ be the elementary pairs where $A$ is moving. Suppose $A=R(x_2,p)\cap R(x_2,q)$. Due to the orientation of $D_0$, one of two things can happen. In case~1, all arcs $B$ of $y$ separating $p,q$ look as on the left in figure~\ref{ga110}.
In case~2, all arcs $B$ of $y$ separating $p,q$ look as on the right. There is always an arc $B$ of $y$ separating $p,q$. Case~1 is equivalent to $(x_3,y)$ being compatible. Case~2 is equivalent to $(x_1,y)$ being compatible. The triple $(x_2,y,A)$ is said to be in {\em case}~1 or {\em case}~2 accordingly.\qed
\end{lemdef}

\newcommand{\hxy}{\pspicture[.5](-9mm, -10mm)(9mm, 9mm) $ \SpecialCoor \psset{unit=.9mm, dimen=middle, radius=1.6pt, labelsep=3pt} \pscircle(0,0){10} 
\pnode(10;70){q0} \pnode(10;110){q1} \pnode(10;-110){q2} \pnode(10;-70){q3}
\Cnode*(-5,0){a} \nput{-90}{a}{p} \Cnode*( 5,0){a} \nput{-90}{a}{q}
$ \endpspicture}

\begin{figure}
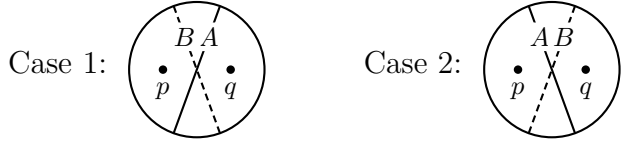

\[ \text{\footnotesize \psset{npos=.24, dash=2.8pt 1.8pt, framesep=2pt} {\normalsize Case 1:\ }
\hxy 
\ncline{q0}{q2} \ncput*{$A$}
\ncline[linestyle=dashed]{q1}{q3} \ncput*{$B$}
\hspace{12mm} {\normalsize Case 2:\ }
\hxy 
\ncline{q1}{q3} \ncput*{$A$}
\ncline[linestyle=dashed]{q0}{q2} \ncput*{$B$}} \]
\caption{\label{ga110} See lemma/definition~\ref{ga59}.}
\end{figure}

\subsection{The universal coverings}

In this intermezzo we introduce some more language. Let $Q\subset\partial D_0$ be an $S$-orbit. We shall define the ``universal coverings'' of our categories $G(Q)$, $G^+(Q)$ written $K(Q)$, $K^+(Q)$.

The object set of $K(Q)$ is $L(Q)$. All hom-sets $K(Q)(x,y)$ consist of one element.

The object set of $K^+(Q)$ is also $L(Q)$. We define $K^+(Q)(x,y)$ to be the set of those elements of $G^+(Mx,My)$ whose images in $G(Mx,My)$ equal $M(x,y)$. The definition of multiplication in $K^+(Q)$ is obvious.
We have a commutative diagram as follows.
\[ \psmatrix[colsep=5.5ex, rowsep=5ex] K^+(Q) & K(Q) \\ G^+(Q) & G(Q) \endpsmatrix
\psset{arrowsize=2pt 4, arrowlength=.7, arrows=->, nodesep=3pt, linewidth=.6pt}
\ncline{1,1}{1,2} \ncline{2,1}{2,2} \ncline{1,1}{2,1} \ncline{1,2}{2,2} \]

Next we look at the presentation for $K^+(Q)$ which is of course closely related to the presentation for $G^+(Q)$.

{\em Generators.} The generators for $K^+(Q)$ are the elementary pairs in $L(Q)\times L(Q)$. The elementary pair $(x,y)$ is a morphism from $x$ to $y$.

A {\em $K$-word\,} (respectively, {\em positive $K$-word}) is a morphism in the groupoid (respectively, category) with object set $L(Q)$ freely generated by the elementary pairs.

{\em Relations.} We call $a_1\cdots a_k=b_1\cdots b_\ell$ a {\em $K(Q)$-relation\,} or $K$-relation (assuming that left and right hand sides are $K(Q)$-words with the same source and target objects) if $(Ma_1)\cdots(Ma_k)=(Mb_1)\cdots (Mb_\ell)$ is an elementary relation (for $G(Q)$). The $K$-words on both sides of the $K$-relation are necessarily positive.

The foregoing taken as a category presentation presents $K^+(Q)$; taken as a groupoid presentation it presents $K(Q)$.

We define an ordering $\leq$ on $L(Q)$ by $x\leq y$ if and only if $K^+(Q)(x,y)$ is non-empty. This is indeed an ordering by \ref{ga61}(b) (atomicity). We haven't yet proved that every hom-set $K^+(Q)(x,y)$ has at most one element or, equivalently, that the natural functor $K^+(Q)\ra K(Q)$ is injective, or again equivalently, that $G^+(Q)\ra G(Q)$ is.

\subsection{Compatibility, part 2}

An elementary pair $(x,x')\in L(Q_1)\times L(Q_1)$ is said to be {\em compatible\,} with $y\in L(Q_2)$ if both $x$ and $x'$ are. Likewise, a $K$-word or $K$-relation is compatible with $y$ if all involved elementary pairs are.

For $y\in L(Q_2)$ let $\leq_y$ denote the ordering on $\Omega(y)$ defined by $x\leq x'$ if and only if there exists a sequence
\[ x=z_0,z_1,\ldots, z_k=x' \]
such that $z_i\in \Omega(y)$ for all $i$ and $(z_{i-1},z_i)$ is an elementary pair. Equivalently, there exists a positive $K(Q_1)$-word from $x$ to $x'$ compatible with~$y$.

It is clear that $x\leq_y x'$ implies $x\leq x'$; the converse will  follow later on in theorem~\ref{ga74}(\ref{ga77}).

\begin{definition} We define a permutation $T$ of $L(Q_1)\cup L(Q_2)$ written on the right. If $\{i,j\}=\{1,2\}$ and $x\in L(Q_i)$, we define $xT\in L(Q_j)$ by rotating the endpoints of all edges of $x$ in positive direction along the boundary of the disk $D_0$ until they first meet $Q_j$.
\end{definition}

Note that for all $y\in L(Q_2)$ one has $yT,yT^{-1}\in\Omega(y)$. These are special elements of $\Omega(y)$ as is shown in the following lemma.

\begin{lemma} \label{ga48} Let $y\in L(Q_2)$. Then $(\Omega(y),\leq_y)$ has a greatest%
\footnote{Let $L$ be an ordered set (not necessarily totally ordered). A {\em greatest element\,} of $L$ is an $x\in L$ such that $y\leq x$ for all $y\in L$. A {\em maximal element\,} of $L$ is an $x\in L$ such that there is no $y\in L$ with $x<y$. Similarly for {\em least\,} and {\em minimal}.}
element $yT$ and a least element $yT^{-1}$.
\end{lemma}

\begin{proof} We shall only be dealing with $\leq_y$, not the ordering on $L(Q)$. Let $x\in \Omega(y)$ be any maximal element. Then for all $A\in\arcs(x)$ the triple $(x,y,A)$ is in case~2 as defined in~\ref{ga59}. A moment's thought now shows that $x=yT$, and we have shown that $\Omega(y)$ has at most one maximal element. Since $\Omega(y)$ is finite by lemma~\ref{ga45}(2) it has a greatest element $yT$. The case of the least element is similar.\end{proof}

\begin{lemma} \label{ga54} Let $y\in L(Q_2)$ and let $a_1\cdots a_k=b_1\cdots b_\ell$ be a $K(Q_1)$-relation. If $a_1$ and $b_1$ are compatible with $y$ then so is the whole relation.
\end{lemma}

\begin{proof} Recall that the $K(Q_1)$-relations are versions of the elementary relations (ER1--4), see figure~\ref{ga26}. First consider the case where our relation corresponds to (ER3), see figure~\ref{ga50}(a). We assume that $y$ is tight with the source object of~$a_1$. The fact that $a_1$ and $b_1$ are compatible with $y$ implies by~\ref{ga59} that some of the arcs of $y$ look like the dashed lines in the figure --- more precisely, their intersections with the three regions that matter. The dashed lines in turn imply, again by~\ref{ga59}, that the whole of the relation is compatible with~$y$.

\newcommand{\lineb}{\psset{linewidth=.8pt, arrowsize=2pt 4, arrowlength=.7, dotstyle=o, dotsize=3.6pt, linestyle=solid, dimen=middle, radius=1.6pt}}
\newcommand{\dashb}{\lineb \psset{linestyle=solid, linewidth=1.8pt, dash=2pt 1.5pt}}
\newcommand{\dashbb}{\psset{linewidth=.8pt, linestyle=dashed, dash=3pt 2pt}}
\newcommand{\ersixdotsb}{{\psset{showpoints=true, linestyle=none} \psline(q0)(q1)(q2)(q3)(q4)(q5)}}
\newcommand{\ersixgonb}{%
{\pspicture(0,0)(0,0) \psset{unit=.08mm} \degrees[120]
\pnode(100; 20){q0}  \pnode(100; 40){q1}
\pnode(100; 64){q2}  \pnode(100; 78){q3}
\pnode(100;102){q4}  \pnode(100;116){q5} 
\psarc(0,0){100}{(q0)}{(q1)} \psarc(0,0){100}{(q2)}{(q3)} \psarc(0,0){100}{(q4)}{(q5)} 
\dashb \psarc(0,0){100}{(q5)}{(q0)} \psarc(0,0){100}{(q1)}{(q2)} \psarc(0,0){100}{(q3)}{(q4)}
\Cnode*(70;8){a} \Cnode*(70;52){a} \Cnode*(55;90){a}
\dashbb \pcarc(100;33)(100;71) \pcarc(100;109)(100;27)
\endpspicture}}
\begin{figure}
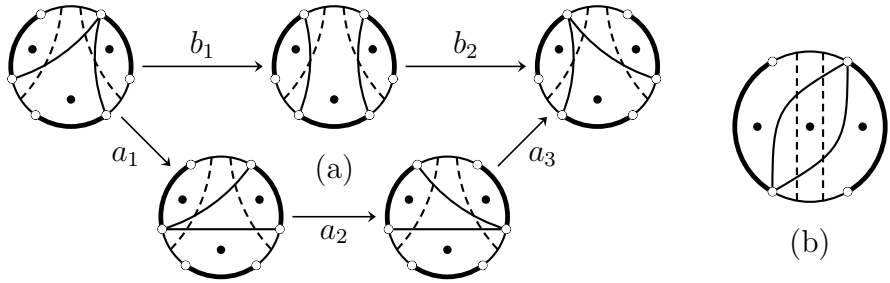
 \lineb \SpecialCoor  \begin{align*}
\psframebox[framesep=8mm, linestyle=none]{\pspicture[.5](0mm,0mm)(70mm,20mm) $ \psset{xunit=3.5mm, yunit=2.5mm} \rput(10,2.5){\text{(a)}}
\pnode( 0,8){n1}  \pnode(10,8){n2}
\pnode(20,8){n3}  \pnode( 5.7,0){n4}
\pnode(14.3,0){n5}
\psset{arcangle=20}
\rput(n1){\ersixgonb} \pcarc(q0)(q2) \pcarc(q4)(q0) \ersixdotsb
\rput(n2){\ersixgonb} \pcarc(q4)(q0) \pcarc(q1)(q3) \ersixdotsb
\rput(n3){\ersixgonb} \pcarc(q1)(q3) \pcarc(q5)(q1) \ersixdotsb
\rput(n4){\ersixgonb} \pcarc(q0)(q2) \psline(q5)(q2) \ersixdotsb
\rput(n5){\ersixgonb} \pcarc(q5)(q1) \psline(q2)(q5)  \ersixdotsb
\psset{nodesep=9.5mm, arrows=->, linewidth=.6pt, labelsep=2.5pt}
\ncline{n1}{n2} \naput{b_1}
\ncline{n2}{n3} \naput{b_2}
\ncline{n5}{n3} \nbput{a_3}
\ncline{n1}{n4} \nbput{a_1}
\ncline{n4}{n5} \nbput{a_2}
$ \endpspicture} \hspace{10mm}
\pspicture[.4](-10mm,-10mm)(10mm,10mm) \psset{unit=.1mm, arcangle=30, ncurv=1.1}
\pnode(100; 60){q0} \pnode(100;120){q1}
\pnode(100;240){q2} \pnode(100;300){q3}
\psarc(0,0){100}{(q0)}{(q1)} \psarc(0,0){100}{(q2)}{(q3)} 
\ncarc{q0}{q2} \ncarc{q2}{q0} 
\dashb \psarc(0,0){100}{(q1)}{(q2)} \psarc(0,0){100}{(q3)}{(q0)}
\psdots(q0)(q1)(q2)(q3)
\Cnode*(0,0){a} \Cnode*(70,0){a} \Cnode*(-70,0){a}
\dashbb \psline(100;100)(100;-100) \psline(100;80)(100;-80)
\uput{10pt}[-90]{0}(100;-90){\text{(b)}}
\endpspicture
\end{align*}
\caption{\label{ga50} See the proof of lemma~\ref{ga54}.}
\end{figure}

In the case (ER4) the fact that $a_1$ and $b_1$ are compatible with $y$ implies that two arcs of $y$ look like the dashed lines in figure~\ref{ga50}(b). We leave it to the reader to conclude that the whole relation is compatible to $y$, as well as to handle the easier cases of (ER1) and (ER2).\end{proof}

\begin{definition} \label{ga88} We define an automorphism $\phi$ of the groupoid $G(Q_1)$. On objects it is defined by $(Mx)\phi=M(xT^2)$. On morphisms it is defined by $(M(x,y))\phi=M(xT^2,yT^2)$.
\end{definition}

We turn to a definition of $\Delta_z$. This is required to be a morphism in $G^+(Q_1)$, not $G(Q_1)$. The attempt $\Delta_{Mx}:=M(x,xT^2)$ would define the correct morphism in $G(Q_1)$, but since we haven't yet proved the functor $G^+(Q_1)\ra G(Q_1)$ to be injective, it wouldn't define $\Delta_z$ as an element of $G^+(Q_1)$. Our next lemma gets around this.

\begin{lemdef} \label{ga57} Let $y\in L(Q_2)$, so that $x:=yT^{-1}$ is the least element of $\Omega(y)$ and $xT^2$ the greatest by lemma~\ref{ga48}. Any positive $K(Q_1)$-word from $x$ to $xT^2$ compatible with $y$ defines the same morphism in $K^+(Q_1)$. The image of this morphism in $G^+(Q_1)$ (from $Mx$ to $Mx\phi$) will be written $\Delta_{Mx}$.
\end{lemdef}

\begin{proof} We only deal with $\leq_y$, not the ordering on $L(Q_1)$. Let $u\in \Omega(y)$ be maximal such that there exist two positive words $a_1\ldots a_k$, $b_1\ldots b_\ell$ compatible with $y$ from $u$ to $xT^2$ defining {\em distinct\,} morphisms in $K^+(Q_1)$. We aim to deduce a contradiction. The maximality assumption implies $a_1\neq b_1$. There is a (unique) $K$-relation $c_1\cdots c_p=d_1\cdots d_q$ with \be c_1=a_1\text{ and }d_1=b_1 \label{ga51} \ee
and we know the whole of this $K$-relation to be compatible with $y$ by lemma~\ref{ga54}.

Extend the positive $K$-word $c_1\cdots c_p$ to a positive $K$-word $c_1\cdots c_r$ ($p\leq r$) compatible with $y$ which ends at $xT^2$. This is possible because $xT^2$ is the greatest element of $\Omega(y)$. Write $\equivp$ for ``defining the same morphism in $K^+(Q_1)$''. Then
\be \begin{aligned}
a_1\cdots a_k &\equivp
(c_1\cdots c_p)(c_{p+1}\cdots c_r) \\ &\equivp
(d_1\cdots d_q)(c_{p+1}\cdots c_r) \equivp
b_1\cdots b_\ell
\end{aligned} \label{ga53} \ee
where the middle $\equivp$ is obvious and the other two equivalences follow from (\ref{ga51}) and maximality of $u$ and the fact that all $K$-words in (\ref{ga53}) are compatible with $y$ by lemma~\ref{ga54}. This contradiction proves the lemma.\end{proof}

\begin{coro} \label{ga68} Let $a\in A(x,y)$ with $x,y\in G_0$. Then there exists $u\in G^+(y,x\phi)$ such that $au=\Delta_x$ in $G^+$.
\end{coro}

\begin{proof}  Write $z=xT$. Note that $y\in\Omega(z)$. Using lemma~\ref{ga48}, there is a positive $K$-word compatible with $z$ from $y$ to $\max\Omega(z)=x\Delta_x$. Its image $u$ in $G^+$ satisfies $au=\Delta_x$.\end{proof}

\begin{definition} Let $Y\subset L(Q_2)$. We define $\Omega(Y):=\cap_{y\in Y}\Omega(y)$ which we endow with an ordering $\leq_Y$ defined by $x\leq_Y x'$ if and only if there exists a positive $K(Q_1)$-word compatible with $Y$ (that is, with all elements of $Y$) from $x$ to $x'$. If $x\leq_Y x'$ then $x\leq_y x'$ for all $y\in Y$.
\end{definition}

\begin{lemma} \label{ga55} Let $(y,z)\in L(Q_2)\times L(Q_2)$ be an elementary pair and put $Y=\{y,z\}$. With respect to $\leq_Y$ then, $\Omega(Y)$ has a greatest element $yT$ and a least element~$zT^{-1}$.
\end{lemma}

\begin{proof} Part (a) of the following figure
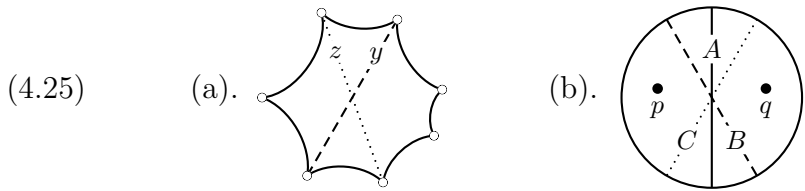
\be \text{(a). \psset{linewidth=.8pt, dotstyle=o, dash=4pt 2.5pt, dotsep=2pt, arcangle=40, labelsep=3pt, unit=.12mm, framesep=2pt} \SpecialCoor
\psframebox[linestyle=none, framesep=0]{\footnotesize 
\pspicture[.5](-13mm,-13mm)(13mm,13mm) $ \psset{npos=.23}
\pnode(100;5){q1} \pnode(100;60){q2} \pnode(100;110){q3} \pnode(100;180){q4} \pnode(100;240){q5} \pnode(100;-70){q6} \pnode(100;-25){q7}
\ncline[linestyle=dashed]{q2}{q5} \ncput*{y} \ncline[linestyle=dotted]{q3}{q6} \ncput*{z}
\pcarc(q1)(q2) \pcarc(q2)(q3) \pcarc[arcangle=50](q3)(q4) \pcarc(q4)(q5) \pcarc(q5)(q6) \pcarc(q6)(q7) \pcarc(q7)(q1)
\psdots(q1)(q2)(q3)(q4)(q5)(q6)(q7)
$\endpspicture} \hspace{10mm} (b).
\psframebox[linestyle=none, framesep=0]{\footnotesize
\pspicture[.5](-13mm,-13mm)(13mm,13mm) $  \psset{radius=2pt, npos=.23}
\pscircle[dimen=middle](0,0){100}
\pcline(100;90)(100;-90) \ncput*{A}
\pcline[linestyle=dashed](100;-60)(100;120) \ncput*{B}
\pcline[linestyle=dotted](100;-120)(100;60) \ncput*{C}
\Cnode*(-60,10){a} \nput{-90}{a}{p}
\Cnode*( 60,10){a} \nput{-90}{a}{q}
$\endpspicture}} \label{ga62} \ee
shows some common edges of $y$ and $z$ (solid) and a dashed arc belonging to $y$ and a dotted arc belonging to $z$.

Let $x\in \Omega(Y)$ be maximal and suppose $x\neq yT$. Since $x\in\Omega(y)$ and $x\neq yT=\max\Omega(y)$ (by lemma~\ref{ga48}) there is an elementary pair $(x,x')$ with $x'\in\Omega(y)$, say obtained by rotating the arc $A$ of~$x$. By \ref{ga59}, some arc $B$ of $y$ looks like the dashed line in (\ref{ga62}b) where $A=R(x,p)\cap R(x,q)$. 

Since $x$ is maximal in $\Omega(Y)$, we have $x'\not\in\Omega(Y)$ and therefore $x'\not\in\Omega(z)$. Also $x\in\Omega(z)$ so by \ref{ga59} some arc $C$ of $z$ looks like the dotted line in (\ref{ga62}b). 

But no pair $\in\arcs(y)\times\arcs(z)$ intersect the way $B$ and $C$ do as dictated by the figure. This is a contradiction. We have proved that $\Omega(Y)$ has at most one maximal element $yT$. Since $\Omega(Y)$ is finite it has a greatest element $yT$. The case of a least element is similar.\end{proof}

\begin{lemma}  \label{ga56} Let $a\in A(x,y)$ be an elementary morphism. Then $a\Delta_y=\Delta_x(a\phi)$.
\end{lemma}

\begin{proof} Let $(u,v)\in L(Q_1)\times L(Q_1)$ be an elementary pair such that $a=M(u,v)$. Then $(uT,vT)\in L(Q_2)\times L(Q_2)$ is an elementary pair too. Put $Y=\{uT,vT\}$. By lemma~\ref{ga55}, $v$ is the least element of $\Omega(Y)$ (relative $\leq_Y$) and $uT^2$ the greatest. 

Choose any positive $K(Q_1)$-word compatible with both elements of $Y$, from $v$ (the least element) to $uT^2$ (the greatest), and let $b$ be the morphism in $K^+(Q_1)$ defined by that word. Let $\pi\col K^+(Q_1)\ra G^+(Q_1)$ be the natural functor. Then $\pi((u,v)b)=\Delta_x$ and $\pi(b(uT^2,vT^2))=\Delta_y$ by \ref{ga57}. Clearly we also have $\pi(u,v)=a$ and $\pi(uT^2,vT^2)=a\phi$. We get
\begin{align*} a\Delta_y
&= \pi(u,v)\cdot \pi\big(b(uT^2,vT^2)\big) \\
&=\pi\big( (u,v)b\big)\cdot \pi(uT^2,vT^2)=
\Delta_x(a\phi). \qedhere \end{align*}\end{proof}

Here is our main theorem.

\begin{theorem} \label{ga66} The tuple $(G,\{A(x,y)\}_{xy},f,\{\Delta_x\}_x,\phi)$ is a Garside tuple.
\end{theorem}

\begin{proof} This follows from \ref{ga61}(c) (atomicity), lemma~\ref{ga70} (cube condition), corollary~\ref{ga68} and lemma~\ref{ga56}.\end{proof}

A {\em lattice\,} is an ordered set $(P,\leq)$ such that any two elements $x,y$ have a least common upper bound (or {\em join}) and a greatest common lower bound (or {\em meet}). 

An ordered set is the same as a category where every hom-set has at most one element, and no two objects are isomorphic. A lattice is then an ordered set in which any two elements have a common upper bound and a common lower bound, and finite limits and colimits exist. This is the language we used in theorem~\ref{ga41}, (\PC\ref{ga35}). It is clear that any two elements of $L$ have a common upper bound and a common lower bound. So a corollary of theorem~\ref{ga66} and (\PC\ref{ga35}) and (\PC\ref{ga37}) is the following.

\begin{coro} \label{ga67} The ordered set $(L,\leq)$ is a lattice.\qed\end{coro}

It is clear from the proof of lemma~\ref{ga6} that Garside's original Garside structure on the braid group \cite{gar69} is a special case of our construction, namely the case where all labels are~2.

\newcommand{\orbitz}{\psdots(p0)(p1)(p2)}

\newcommand{\orbity}{%
\raisebox{.7ex}{\pspicture[.5](-6,-6)(6,6) \pscircle(0,0){5}
\pnode(5;14){p0} \pnode(5;22){p1} \pnode(5;6){p2} \pnode(7.5;3){q}
\endpspicture}}

\newcommand{\orbita}{\orbity \ncline{p1}{p2} \ncarc{p1}{p2} \orbitz}
\newcommand{\orbitb}{\orbity \ncline{p2}{p0} \ncarc{p2}{p0} \orbitz}
\newcommand{\orbitc}{\orbity \psline(p0)(p2)(p1) \orbitz}

\begin{example} Let $x\in L$. If the orbit $Mx$ is a sublattice of $L$ (that is, closed under join and meet) we obtain a lattice ordering on $M$ by putting $g\leq h$ $\Leftrightarrow$ $gx\leq hx$. But $M$-orbits in $L$ are not sublattices of $L$ in general, as is shown by the following example.
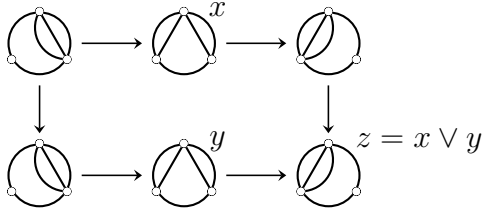
\[ \SpecialCoor \psset{unit=.85mm, linewidth=.8pt, dotstyle=o, dotsize=3.4pt, arcangle=60, ncurv=.9} \degrees[24]
\psmatrix[colsep=9mm, rowsep=7mm]
\orbita & \orbitc \rput(q){x} & \orbitb \\
\orbita & \orbitc \rput(q){y} & \orbitb  \rput(q){\makebox[1ex][l]{$z=x\vee y$}}
\endpsmatrix
\psset{arrowsize=2pt 4, arrowlength=.7, arrows=->, nodesep=.5mm, linewidth=.6pt}
\ncline{1,1}{1,2} \ncline{1,2}{1,3}
\ncline{2,1}{2,2} \ncline{2,2}{2,3}
\ncline{1,1}{2,1} \ncline{1,3}{2,3} 
\]
This diagram is an elementary relation, so $z=x\vee y$. Now $x,y$ are in the same orbit but $x\vee y$ is in a different orbit. 
\end{example}

\section{Finite lattices and associahedra} \label{ga92}

\newcommand{\dashd}{\psset{linewidth=1.8pt}}

\newcommand{\sixdotsd}{\psdots(q0)(q6)(q12)(q18)(q24)(q30)}

\newcommand{\sixgond}{%
{\pspicture[.4](-7mm,-7mm)(7mm,7mm) \psset{unit=7mm, dimen=middle, linestyle=solid} \degrees[36]
\pnode(1;0){q0} \pnode(1;3){q3} \pnode(1;6){q6} \pnode(1;8){q8} \pnode(1;10){q10} \pnode(1;12){q12} \pnode(1;15){q15} \pnode(1;18){q18} \pnode(1;20){q20} \pnode(1;22){q22} \pnode(1;24){q24} \pnode(1;30){q30} \pnode(1;32){q32} \pnode(1;34){q34}
\psarc(0,0){1}{(q0)}{(q6)} \psarc(0,0){1}{(q12)}{(q18)} \psarc(0,0){1}{(q24)}{(q30)} 
\dashd \psarc(0,0){1}{(q6)}{(q12)} \psarc(0,0){1}{(q18)}{(q24)} \psarc(0,0){1}{(q30)}{(q0)}
\endpspicture}}

Recall that in corollary~\ref{ga67} we constructed a lattice $L(Q_1)$. The braid-like group $M$ acts freely on it with finitely many orbits. We also have the finite subsets $\Omega(y)\subset L(Q_1)$ where $y\in L(Q_2)$. If all labels are 2, then $L(Q_1)$ can be identified with the braid group $B_n$ and $\Omega(y)$ with the symmetric group $S_n$.

In theorem~\ref{ga74} we prove that $\Omega(y)$ is a sublattice of $L(Q_1)$. In particular, it is a lattice in its own right. Most particular cases of these finite lattices seem to be new.

It seems that $\Omega(y)$ is the vertex set of a natural polytope (up to deformation); we shall not go into this. If all labels are 2 it is known as the permutahedron.

In subsection~\ref{ga87} we study the case where all labels are $3$ more closely.

The {\em associahedron\,} is a certain polytope whose vertex set is the set of triangulations of a fixed $n$-gon such that the vertex set of the triangulation is the vertex set of the $n$-gon \cite{sta63}, \cite{lee89}. 

Suppose now that all labels are~3. Recall that the object set of the groupoid $G=G(Q_1)$ is $G_0=M\backslash L(Q_1)$. This is essentially the set of triangulations of the disk $D_0$ with vertex set $Q_1$, that is, the vertex set of the associahedron.

In proposition~\ref{ga79} we prove that $\Omega(y)$ is in bijection with $G_0$. Combined with the lattice ordering on $\Omega(y)$ mentioned above this yields a family of lattice orderings on the vertex set of the associahedron. One of these orderings is known as the Tamari lattice.

We have already come across pictures of our orderings on the  vertex set of the $3$-dimensional associahedron in figure \ref{ga42}~(b)--(d) when we were studying characteristic graphs.  The 3-dimensional Tamari lattice is part~(c) of the figure.

In contrast, we obtain essentially only one lattice ordering on the vertex set of the permutahedron, because if all labels are 2 then the $M$-action on $L$ is transitive.

\subsection{More on $\Omega(y)$}

\begin{theorem} \label{ga74} Let $y\in L(Q_2)$.%
\lista{{\rm(\arabic{daana})}}{\setlength{\leftmargin}{3em}}
\item \label{ga75} Let $a_1\cdots a_k=b_1\cdots b_\ell$ be a $K(Q_1)$-relation {\rm (}in particular, the $a_i$ and $b_j$ are elementary pairs in $L(Q_1)\times L(Q_1))$. Then $a_1\cdots a_k$ is compatible with $y$ if and only if $b_1\cdots b_\ell$ is.
\item \label{ga76} We have $\Omega(y)=\{x\in L(Q_1): yT^{-1}\leq x\leq yT\}$. 
\item \label{ga77} Let $x,x'\in\Omega(y)$. Then $x\leq_y x'$ is equivalent to $x\leq x'$.
\item \label{ga78} The ordered set $(\Omega(y),\leq)$ is a lattice.
\end{list}
\end{theorem}

\begin{proof} (1). Suppose the relation is of the form (ER3) and more precisely reads $a_1a_2a_3=b_1b_2$ as in figure~\ref{ga26}. We will show that if $a_1a_2a_3$ is compatible with $y$ then so is $b_1b_2$. We apply \ref{ga59} three times as follows. Since $a_1$ is compatible with $y$, at least one of the dashed lines in part (a) of the following figure
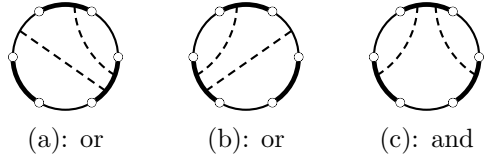
\[ \text{\footnotesize \psset{dotstyle=o, dotsize=3.6pt, linewidth=.8pt, arcangle=30, linestyle=dashed, dash=3pt 2pt} \SpecialCoor
\psmatrix[colsep=10mm, rowsep=2mm]
\sixgond \ncline{q15}{q32} \ncarc{q34}{q8} \sixdotsd &
\sixgond \ncline{q22}{q3} \ncarc{q10}{q20} \sixdotsd &
\sixgond \ncarc{q34}{q8} \ncarc{q10}{q20} \sixdotsd \\
(a): or & (b): or & (c): and
\endpsmatrix} \]
describes an arc in~$y$. Since $a_3$ is compatible with $y$, at least one of the dashed lines in part (b) also describes an arc in~$y$. Since arcs of $y$ don't intersect in the interior of the disk $D_0$ {\em both\,} dashed lines in part (c) describe arcs of~$y$. They in turn imply that $b_1b_2$ is compatible with~$y$. They also imply that figure~\ref{ga50}(a) matches our situation, including the dashed lines.

The reverse implication for (ER3) as well as the remaining elementary relations are similar and left to the reader. Note that (ER4) is slightly harder than (ER3) because the punctures cannot be omitted.

(2). The inclusion $\subset$ follows from lemma~\ref{ga55}. We shall prove $\supset$. Let $X$ be the set of positive $K$-words from $yT^{-1}$ to $yT$. Let $d$ be the largest metric (in a generalised sense, namely distances can be $\infty$) such that $d(uv_1w,uv_2w)=1$ if ``$v_1=v_2$'' is a $K$-relation. In fact, any two elements of $X$ have finite distance because the functor $G^+\ra G$ is injective by (\PC\ref{ga37}) and our main theorem~\ref{ga66}. But (1) states that if two elements of $X$ have distance~$1$ and one is compatible with $y$ then so is the other. Since at least one element of $X$ is compatible with $y$ (by lemma~\ref{ga55}), all are as required.

(3). The implication $\Rightarrow$ is trivial. We prove $\Leftarrow$. Let $x\leq x'$. Then there exists a positive $K$-word $w$ from $x$ to $x'$. By part (2), $w$ is compatible with $y$. It follows that $x\leq_y x'$.

(4). This follows directly from (2) and the fact that $(L(Q_1),{\leq})$ is a lattice by corollary~\ref{ga67}.\end{proof}

\subsection{The case of only triangles} \label{ga87}

In this subsection, we assume all labels to be~3. 

\begin{prop} \label{ga79}  Suppose all labels are~$3$. For any $y\in L(Q_2)$, the natural map $\pi\col\Omega(y)\ra G_0$ is bijective.
\end{prop}

\begin{proof} Injectivity is easy and left to the reader (and true even if all labels are  $\geq 3$).

We shall prove surjectivity. Let $x'\in G_0$, $y\in L(Q_2)$ be tight. 

Let $R$ be a region of $y$ with vertices $r_1,r_2,r_3$. Among all arcs of $x'$ intersecting $R$ and separating $r_i$ from the other two vertices ($i=1,2,3$)%
\footnote{There is at least one such arc because of the boundary edges of $x'$.},
let $A_i$ be the innermost one. By $R'$ we shall denote the region of $x'$ containing $A_1,A_2,A_3$. 

There exists a self-homeomorphism of $D_0$ fixing the boundary pointwise which preserves $y$ and which takes the puncture in any region $R$ of $y$ to an interior point of $R\cap R'$. Then $x:=h^{-1}x'$ is in $\Omega(y)$ (using the fact that all labels are $3$); it is clear that $\pi x=x'$. This proves surjectivity.\end{proof}

Define $\pi$ as in proposition~\ref{ga79}. For any $x\in L(Q_1)$ we define an ordering $\leq_{\pi x}$ on $G_0$ as follows. For $y,z\in\Omega(xT)$ we put \[ \pi y\leq_{\pi x} \pi z\Longleftrightarrow y\leq z \]
which is equivalent to $y\leq_{xT} z$ by theorem~\ref{ga74}(\ref{ga77}). Clearly, $\leq_{\pi x}$ depends only on $\pi x$. 

Notice that $x\leq_x y$ for all $x,y\in G_0$, and that $\leq_x$ is a lattice ordering on $G_0$, the vertex set of the associahedron. So we have many lattice orderings $\leq_x$ on $G_0$. If $x$ is a fan, this ordering was discovered by Tamari \cite{tam51}, \cite{fritam67}, \cite[page~18]{gra78}. If $x$ is another triangulation of the $n$-gon it seems to be new.

We get the following amusing property which we shall not prove. The braid-like groupoid $G$ with object set $G_0$ is presented by a generator $[xy]\in G(x,y)$ (for all $x,y\in G_0$) and relations $[xx]=1$ for all $x$ and $[xy][yz]=[xz]$ whenever $y\leq_x z$.


\begin{thebibliography}{B}

\bibitem[BKL98]{bkl98} Birman, Joan;  Ko, Ki Hyoung;  Lee, Sang Jin. A new approach to the word and conjugacy problems in the braid  groups. Adv.\ Math.\  {\bf 139}  (1998),  no. 2, 322--353.


\bibitem[Deh00]{deh00} Dehornoy, Patrick. Chapter~2 in {\em Braids and self-distributivity}.
Progress in Mathematics, 192. Birkh\"auser Verlag, Basel,  2000.

\bibitem[Deh02]{deh02} Dehornoy, Patrick. Groupes de Garside.  Ann.\ Sci.\ \'Ecole Norm.\ Sup.\ (4)  {\bf 35}  (2002),  no.\ 2, 267--306.

\bibitem[DehPar99]{dehpar99} Dehornoy, Patrick;  Paris, Luis. Gaussian groups and Garside groups, two generalisations of Artin  groups. Proc.\ London Math.\ Soc.\ (3)  {\bf 79}  (1999),  no. 3, 569--604.

\bibitem[Eps92]{eps92} Epstein, D.B.A.; J.W.\ Cannon; D.F.\ Holt; S.V.F.\ Levy; M.S.\ Paterson; W.P.\ Thurston. Word processing in groups. Jones and Bartlett Publishers, Boston, MA, 1992.

\bibitem[FriTam67]{fritam67} Friedman, Haya;  Tamari, Dov. Probl\`emes d'associativit\'e: Une structure de treillis finis induite  par une loi demi-associative. J.\ Combinatorial Theory  {\bf 2}  (1967), 215--242.

\bibitem[Gar69]{gar69} Garside, F.\ A.  The braid group and other groups. Quart.\ J.\ Math.\ Oxford Ser.\ (2)  {\bf 20}  (1969), 235--254.

\bibitem[Gr\"a78]{gra78}Gr\"atzer, George. {\em General lattice theory}. Second edition. Birkh\"auser Verlag, Basel, 1998. (First edition published 1978).

\bibitem[Lee89]{lee89} Lee, Carl W.  The associahedron and triangulations of the $n$-gon. European J.\ Combin.\  {\bf 10}  (1989),  no.\ 6, 551--560.

\bibitem[Par05]{par05} Paris, Luis. From braid groups to mapping class groups. \hfill\\ {\tt http://arxiv.org/abs/math.GR/0412024}.

\bibitem[Sta63]{sta63} Stasheff, James. Homotopy associativity of $H$-spaces I. Trans.\ Amer.\ Math.\ Soc.\ {\bf 108} (1963), 275--292.

\bibitem[Tam51]{tam51}
Tamari, Dov. Mono\"ides pr\'eordonn\'es et chaines de Malcev. Thesis, Paris, 1951.

\end{thebibliography}
\end{document}